\documentclass[final,leqno]{siamltex}
\usepackage{amssymb}
\usepackage{amsmath}
\usepackage{color}
\usepackage{enumerate}
\usepackage{pgfplots}
\usepackage{tikz}
\usetikzlibrary{shapes.geometric}

\usepackage{graphicx}

\newcommand{\cT}{\mathcal{T}}

\newcommand{\cM}{\mathcal{M}}

\newcommand{\NJ}{\mathbb{N}_J^d}

\newcommand{\bx}{\mathbf{x}}

\newcommand{\del}{\delta}

\newcommand{\bh}{\mathbf{h}}

\newtheorem{remark}{Remark}[section]

\allowdisplaybreaks[1]

\begin{document}

\title{A high order correction to the Lax-Friedrich's method for approximating 
stationary Hamilton-Jacobi equations}
 
\author{
Thomas Lewis\thanks{Department of Mathematics and Statistics, 
The University of North Carolina at Greensboro, 
Greensboro, NC 27402, U.S.A.  {\tt tllewis3@uncg.edu}.
The work of this author was partially supported by the NSF grant DMS-2111059.}
\and
Xiaohuan Xue\thanks{Department of Mathematics and Statistics, 
The University of North Carolina at Greensboro, 
Greensboro, NC 27402, U.S.A.  {\tt x\_xue2@uncg.edu}.
The work of this author was partially supported by the NSF grant DMS-2111059.}
}

\maketitle

\begin{abstract} 
A new class of non-monotone finite difference (FD) approximation methods 
for approximating solutions to non-degenerate stationary Hamilton-Jacobi 
problems with Dirichlet boundary conditions 
is proposed and analyzed.
The new FD methods add a high order correction to the Lax-Friedrich's method while 
utilizing a novel cutoff to preserve the convergence properties of the Lax-Friedrich's approximation.  
Since monotone methods are limited to first order accuracy by the Godunov barrier, 
the proposed approach provides a template for boosting the accuracy of a monotone method 
using a modified numerical moment stabilizer with a high-order auxiliary boundary condition.  
Numerical tests are provided to test the utility of the approach while a novel 
admissibility and stability analysis technique lays a foundation for 
analyzing non-monotone methods.  
\end{abstract}



\pagestyle{myheadings}
\thispagestyle{plain}
\markboth{T. Lewis and X. Xue}{A high order correction of the Lax-Friedrich's method}


\section{Introduction}
In this paper we propose two new finite difference (FD) methods for approximating 
the viscosity solution of the stationary Hamilton-Jacobi (HJ) equation 
\begin{subequations} \label{HJ}
\begin{alignat}{2}
	H[u] \equiv H(\nabla u , u , \bx) + \theta u & = 0 && \text{in } \Omega , \label{HJ_pde} \\ 
	u & = g \qquad && \text{on } \partial \Omega , \label{HJ_bc}
\end{alignat}
\end{subequations}
where 
$\theta > 0$ is a constant; 
$\nabla u (\bx)$ is the gradient of $u$ at $\bx$;  
$\Omega \subset \mathbb{R}^{d}$ is an open, bounded domain; 
$g$ is continuous on $\partial \Omega$; 
and $H$ is locally Lipschitz and nondecreasing with respect to $u$. 
Under these assumptions, \eqref{HJ} has a unique viscosity solution in the space 
$C(\overline{\Omega})$. 
The new, non-monotone FD methods will be capable of exceeding the first-order 
accuracy barrier associated with monotone FD methods by adding a 
correction term to the first-order accurate Lax-Friedrich's method.  
Numerical tests show that the methods are capable of 
approximating classical solutions of \eqref{HJ} with up to second order accuracy as well as 
reliably approximating lower-regularity viscosity solutions of \eqref{HJ}.  

The formulations and analysis in this paper can trivially be extended to 
the vanishing viscosity regularization of \eqref{HJ} defined by
\begin{subequations}\label{HJ_vis}
\begin{alignat}{2}
    H^\epsilon[u^\epsilon] \equiv - \epsilon \Delta u^\epsilon + H(\nabla u^\epsilon, u^\epsilon , \bx) 
    	+ \theta u^\epsilon & = 0 && \text{in } \Omega , \\ 
	u^\epsilon & = g \qquad && \text{on } \partial \Omega 
\end{alignat}
\end{subequations}
for $\epsilon > 0$, which is a proper elliptic second order problem. 
With appropriate assumptions for $H^\epsilon$, the solution $u^\epsilon$ exists 
and $u^\epsilon$ converges to the viscosity solution of \eqref{HJ} in $L^\infty (\Omega)$ 
at a rate of $\mathcal{O}(\sqrt{\epsilon})$. 
The function $u^\epsilon$ is often called the vanishing viscosity approximation of the viscosity solution 
of \eqref{HJ} for small $\epsilon$. 
See Section~\ref{visc_sec} for a more detailed discussion regarding 
viscosity solutions and the existence and uniqueness theory for \eqref{HJ} and \eqref{HJ_vis} 
as well as how the Dirichlet boundary condition can be interpreted in the viscosity sense 
without specifying inflow and outflow regions for the boundary. 

The schemes proposed in this paper use a numerical moment stabilizer 
to increase their accuracy when compared to a vanishing viscosity approach.  
In order to prove admissibility, stability, and convergence 
properties of an underlying monotone approximation scheme can be preserved, 
a cutoff operator is introduced to account for the non-monotone nature of the numerical 
moment.  
The main idea is to add correction terms that appear in the numerical moment 
to the Lax-Friedrich's method as a way to decrease the 
numerical diffusivity error.  
Unfortunately, standard analytic techniques for approximating solutions to fully nonlinear degenerate elliptic 
problems require monotonicity.    
The cutoff for the numerical moment is introduced to account for the non-monotone structure 
of the proposed scheme along with a modified sub- and supersolution approach for proving 
admissibility and uniform stability results.  
Future work will seek to remove the need for the cutoff in the analysis; however, the cutoff can be tuned 
to yield approximations consistent with the proposed second-order 
scheme that uses the unmodified numerical moment stabilizer.  

Stationary Hamilton-Jacobi equations arise from various scientific applications including 
optimal control, wave propagation, geometric optics, multiphase flow, image processing, 
etc. (cf. \cite{Osher, Sethian99} and the references therein). 
The numerical approximations to the solutions of time-dependent and stationary Hamilton-Jacobi equations 
have been significantly studied 
(cf. \cite{Barth, Crandall_Lions84, CockburnD, Hu, Huang, Kao, OsherShu, Sethian03, Shu, Tadmor, Tsai} 
and the references therein) 
as they are vital in understanding application problems. 
Monotone FD methods under the framework of Crandall and Lions in \cite{Crandall_Lions84} 
like the Lax-Friedrich's method borrowed from the approximation theory for nonlinear hyperbolic conservation laws 
(see \cite{Tadmor87}) have formed an analytic foundation for approximating viscosity solutions. 
Stationary Hamilton-Jacobi equations can be viewed  as degenerate fully nonlinear elliptic problems. 
The theory for monotone FD methods was further extended to elliptic problems by Barles and Souganidis in \cite{Barles}, 
which guarantees convergence to the underlying viscosity solution 
if the PDE satisfies the comparison principle and if the numerical scheme is 
monotone, admissible, consistent, and stable. 
Unfortunately, monotone methods for approximating \eqref{HJ} 
are limited to first-order accuracy due to the Godunov barrier (see \cite{Tadmor}). 
This paper proposes a new, non-monotone numerical method for problem \eqref{HJ} with potential 
to achieve higher-order accuracy. 
Note that several formal approaches have been used to design high-order methods to overcome the Godunov barrier
as seen in \cite{Barth, Shu, centralDG_Li,Xu,Yan_Osher11} and the references therein;  
however, convergence results for higher order methods typically are open due to the lack of monotonicity 
or hold for dynamic problems but do not extend to stationary problems.  
The primary hurdle for the stationary case is developing analytic techniques for proving stability 
without assuming the scheme is monotone.  
This paper introduces a novel cutoff technique and analytic technique for 
formulating and analyzing a class of high-order, non-monotone, stable schemes 
in the stationary case with applications to the dynamic case.  
The cutoff technique ensures the non-monotone terms that increase the order of accuracy 
do not push the approximation too far from an underlying converging monotone method 
that reliably approximates the viscosity solution.  

Similar to the time-dependent Hamilton-Jacobi equations discretized by appropriate time-stepping schemes, 
the stationary problem \eqref{HJ} has some inherent difficulties in establishing the admissibility and stability analysis.  
In contrast, the regularized problem \eqref{HJ_vis} is uniformly elliptic, 
and the extra viscous term for a fixed positive $\epsilon$ can help to control advective terms 
when the mesh size is sufficiently small.   
When viewing \eqref{HJ_vis} as a viscosity approximation to \eqref{HJ}, 
requiring the mesh to resolve the viscosity term in order to ensure monotonicity 
presents a strong mesh condition and/or creates additional error based on the scale of $\epsilon$.  
The techniques in this paper provide a new approach for  
analyzing the admissibility, stability, and convergence of 
a non-monotone method for approximating \eqref{HJ_vis} 
without introducing a mesh restriction 
and without being limited by a direct realization of the vanishing viscosity approach.  
The key concept in the admissibility and stability analysis is to recast the method as the solution to 
a fixed-point problem and employ the Schauder fixed-point 
theorem instead of the stronger contractive mapping theorem that is typically employed 
when analyzing monotone methods.  
Thus, we can weaken the sufficient conditions that ensure the existence of a fixed point.  
The fixed-point will be bounded above and below by a strategic choice of sub- and supersolutions 
formed using a monotone approximation method.  

The main idea in this paper will be to treat \eqref{HJ} as a degenerate second order operator and  
introduce a {\em numerical moment} that will serve 
as a higher order stabilizer, as seen in \cite{FDhjb} and \cite{Kellie} for fully nonlinear uniformly elliptic 
second-order problems 
and \cite{FD_CDR} for linear, constant-coefficient first- and second-order problems.  
Using a numerical moment with appropriate auxiliary boundary conditions 
instead of a numerical viscosity (as typical for monotone methods) 
will allow us to break the first-order accuracy barrier. 
Our proposed method allows for local truncation errors to achieve up to second-order accuracy 
thanks to the introduction of 
a high-order auxiliary boundary condition and appropriate scaling of the numerical moment 
that was not considered in \cite{FDhjb} and \cite{Kellie} 
and only numerically tested in \cite{FD_CDR}.  
Unfortunately, the analytic techniques for the admissibility and stability analysis in \cite{FDhjb} and \cite{Kellie} 
do not directly translate to the method proposed for approximating \eqref{HJ} 
since they did not account for the presence of the gradient operator as observed in \cite{FD_CDR}.  
Thus, an appropriate cutoff technique is introduced.  
A benefit of the cutoff is that it will directly yield $\ell^\infty$ stability bounds 
as opposed to the discrete Sobelev inequality used in \cite{FDhjb} and \cite{Kellie} that 
required the number of spatial dimensions to be limited to $d \leq 3$.  

The remainder of this paper is organized as follows. 
In Section~\ref{visc_sec} we introduce some preliminaries which include background for viscosity solutions, 
Hamilton-Jacobi equations, and monotone methods as well as notation. 
We formulate the new FD methods in Section~\ref{formulation_sec}.  
We first introduce the scheme that uses a numerical moment while proving various 
consistency results and providing motivation for the formulation.  
We then introduce the scheme that uses a cutoff to modify the numerical moment 
while still providing a potentially high-order correction to the Lax-Friedrich's method.  
We discuss various admissibility, stability, and convergence results for 
the Lax-Friedrich's method, the proposed high-order method, and the modified version with a cutoff 
in Section~\ref{analysis_sec}.
Several numerical tests are provided in Section~\ref{numerics_sec} to verify the accuracy
of the proposed methods, 
and some concluding remarks can be found in Section~\ref{conclusion_sec}.


\section{Preliminaries} \label{visc_sec}

The appropriate solution concept for \eqref{HJ} is viscosity solution theory.  
We begin by introducing some basics for viscosity solution theory 
such as the notation, the definition of a solution, and the comparison principle.  
We next introduce the difference operator notation that will be used to formulate 
the approximation schemes.    
We end the section by recalling the Lax-Friedrich's method and some standard results 
for monotone methods.  


\subsection{Viscosity solutions} \label{spaces_sec}

Problem \eqref{HJ} is not guaranteed to have a classical solution in 
$C^1(\Omega) \cap C (\overline{\Omega})$ 
and problem \eqref{HJ_vis} is not guaranteed to have a classical solution in 
$C^2(\Omega) \cap C (\overline{\Omega})$.  
Furthermore, there may exist multiple functions that satisfy the PDE almost everywhere in $\Omega$.  
The primary strength of viscosity solution theory is that it provides a framework weak enough 
to allow solutions of \eqref{HJ} in $C (\overline{\Omega})$ 
while also having enough structure to guarantee the uniqueness of the solution in many instances.  
A secondary strength is that it allows the Dirichlet condition to be imposed over the 
entire boundary since the inflow and outflow boundary would depend on the unknown solution $u$. 
Since we are assuming both $H$ and $g$ are continuous in \eqref{HJ}, 
we have the following definition of a viscosity solution 
that satisfies the boundary condition in the viscosity sense:   

\begin{definition} \label{visc_def}
Let $H$ denote the differential operator in \eqref{HJ_pde} 
and $g$ denote the function in \eqref{HJ_bc}.  
\begin{enumerate}[{\rm (i)}]
\item
A function $u \in C(\overline{\Omega})$ 
is called a {\em viscosity subsolution} of \eqref{HJ}
if $\; \forall \varphi\in C^1 (\overline{\Omega})$, when $u-\varphi$ has a local 
maximum at $\bx_0 \in \overline{\Omega}$ with $u(\bx_0) = \varphi(\bx_0)$, 
\[
H(\nabla \varphi(\bx_0) , \varphi(\bx_0) , \bx_0) \leq 0 
\]
if $\bx_0 \in \Omega$ or 
\[
\varphi(\bx_0) - g(\bx_0) \leq 0 
\qquad \text{or} \qquad 
H(\nabla \varphi(\bx_0) , \varphi(\bx_0) , \bx_0) \leq 0 
\]
if $\bx_0 \in \partial \Omega$.  
\item 
A function $u \in C(\overline{\Omega})$ 
is called a {\em viscosity supersolution} of \eqref{HJ}
if $\; \forall \varphi\in C^1 (\overline{\Omega})$, when $u-\varphi$ has a local 
minimum at $\bx_0 \in \Omega$ with $u(\bx_0) = \varphi(\bx_0)$, 
\[
H(\nabla \varphi(\bx_0) , \varphi(\bx_0) , \bx_0) \geq 0 
\]
if $\bx_0 \in \Omega$ or 
\[
\varphi(\bx_0) - g(\bx_0) \geq 0 
\qquad \text{or} \qquad 
H(\nabla \varphi(\bx_0) , \varphi(\bx_0) , \bx_0) \geq 0 
\]
if $\bx_0 \in \partial \Omega$.  
\item 
A function $u \in C(\overline{\Omega})$ 
is called a {\em viscosity solution} of \eqref{HJ}
if $u$ is both a viscosity subsolution and a viscosity supersolution of \eqref{HJ}.
\end{enumerate}
\end{definition}

A constructive way to understand viscosity solutions of \eqref{HJ} is as the limit 
of the solutions to a family of quasilinear second order problems (\cite{Crandall_Lions84}).  
Consider the regularized problem \eqref{HJ_vis}.    
The solution $u^\epsilon$ converges uniformly to the solution $u$ of \eqref{HJ}  
with the term $- \epsilon \Delta u^\epsilon$ referred to as a vanishing viscosity.  

The nonlinear differential operator $H$ in \eqref{HJ_pde} is a degenerate elliptic operator 
in the sense that, trivially, it is nonincreasing with respect to the Hessian $D^2 u$ 
using the natural partial ordering of $\mathcal{S}^{d \times d}$, the space of symmetric matrices, 
where $A \geq B$ if $A-B$ is a nonnegative definite matrix.  
Furthermore, problem \eqref{HJ} with $H$ and $g$ continuous and $\theta > 0$ has a unique solution 
since it satisfies a comparison principle 
as stated in Definition~\ref{comparison}.  

\begin{definition}\label{comparison}
Problem \eqref{HJ} is said to satisfy a {\em comparison 
principle} if the following statement holds. For any upper semi-continuous 
function $u$ and lower semi-continuous function $v$ on $\overline{\Omega}$,  
if $u$ is a viscosity subsolution and $v$ is a viscosity supersolution 
of \eqref{HJ}, then $u\leq v$ on $\overline{\Omega}$.
\end{definition} 


\subsection{Difference operators} 

We introduce various difference operators for approximating first and second order partial derivatives.  
Assume $\Omega$ is a $d$-rectangle, i.e., 
$\Omega = \left( a_1 , b_1 \right) \times \left( a_2 , b_2 \right) \times \cdots \times 
		\left( a_d , b_d \right)$.    
We will only consider grids that are uniform in each coordinate $x_i$, $i = 1, 2, \ldots, d$, 
in this paper; however, the formulations in Section~\ref{formulation_sec} and analytic techniques in Section~\ref{analysis_sec} 
can be extended to non-uniform grids.  
Let $J_i$ be a positive integer and $h_i = \frac{b_i-a_i}{J_i-1}$ for $i = 1, 2, \ldots, d$. 
Define $\mathbf{h} = \left( h_1, h_2, \ldots, h_d \right) \in \mathbb{R}^d$, 
$h = \max_{i=1,2,\ldots,d} h_i$, 
$J = \prod_{i=1}^d J_i$, and   
$\NJ = \{ \alpha = (\alpha_1, \alpha_2, \ldots, \alpha_d) 
\mid 1 \leq \alpha_i \leq J_i, i = 1, 2, \ldots, d \}$.  Then, $\left| \NJ \right| = J$. 
We partition $\Omega$ into $\prod_{i=1}^d \left(J_i-1 \right)$ sub-$d$-rectangles with grid points
$\bx_{\alpha} = \Big(a_1+ (\alpha_1-1)h_1 , a_2 + (\alpha_2-1) h_2 , \ldots , 
		a_d + (\alpha_d-1) h_d \Big)$
for each multi-index $\alpha \in \NJ$.
We call $\cT_{\mathbf{h}}=\{ \bx_{\alpha} \}_{\alpha \in \NJ}$ 
a grid (set of nodes) for $\overline{\Omega}$. 

Let $\left\{ \mathbf{e}_i \right\}_{i=1}^d$ denote the canonical basis vectors for $\mathbb{R}^d$. 
We also introduce an interior grid that removes the layer of grid points adjacent to the boundary of $\Omega$.  
To this end, let $\mathring{\cT_{\bh}} \subset \cT_{\bh}$ such that 
\[
	\mathring{\cT_{\bh}} \equiv 
	\left\{ \bx_\alpha \in \cT_{\bh} \cap \Omega \mid \bx_{\alpha \pm \mathbf{e}_i} \notin \partial \Omega 
	\text{ for all } i = 1,2,\ldots,d \right\} . 
\]

Define the (first order) forward and backward difference operators by
\begin{equation} \label{fd_x}
	\del_{x_i,h_i}^+ v(\bx)\equiv \frac{v(\bx + h_i \mathbf{e}_i) - v(\bx)}{h_i},\qquad
	\del_{x_i,h_i}^- v(\bx)\equiv \frac{v(\bx)- v(\bx-h_i \mathbf{e}_i)}{h_i}
\end{equation}
for a function $v$ defined on $\mathbb{R}^d$ and 
\[
	\del_{x_i,h_i}^+ V_\alpha \equiv \frac{V_{\alpha + \mathbf{e}_i} - V_\alpha}{h_i} , \qquad
	\del_{x_i,h_i}^- V_\alpha \equiv \frac{V_\alpha- V_{\alpha - \mathbf{e}_i}}{h_i}
\]
for a grid function $V$ defined on the grid $\mathcal{T}_{\mathbf{h}}$.  
We also define the following (second order) central difference operator: 
\begin{equation} \label{fd_xc}
	\delta_{x_i, h_i}  \equiv \frac{1}{2} \left( \delta_{x_i, h_i}^+ + \delta_{x_i, h_i}^- \right) 
\qquad\mbox{for } i=1,2,\cdots, d 
\end{equation}
and the ``sided" and central gradient operators $\nabla_{\mathbf{h}}^+$, 
$\nabla_{\mathbf{h}}^-$, and $\nabla_{\mathbf{h}}$ by 
\begin{align} \label{discrete_grad_def}
	\nabla_{\mathbf{h}}^\pm \equiv \bigl[ \delta^\pm_{x_1, h_1} , \delta_{x_2, h_2}^\pm , \cdots , 
		\delta_{x_d, h_d}^\pm \bigr]^T, \quad   
	\nabla_{\mathbf{h}} \equiv \bigl[ \delta_{x_1, h_1} , \delta_{x_2, h_2} , \cdots , 
		\delta_{x_d, h_d} \bigr]^T. 
\end{align}

Define the (second order) central difference operator for approximating second order 
non-mixed partial derivatives by 
\begin{equation} \label{fd_xx}
	\del_{x_i,h_i}^2 v(\bx) 
	\equiv \frac{v(\bx + h_i \mathbf{e}_i) - 2 v(\bx) + v(\bx - h_i \mathbf{e}_i)}{h_i^2} 
\end{equation}
for a function $v$ defined on $\mathbb{R}^d$ and 
\[
	\del_{x_i,h_i}^2 V_\alpha 
	\equiv \frac{V_{\alpha + \mathbf{e}_i} - 2 V_\alpha + V_{\alpha - \mathbf{e}_i}}{h_i^2} 
\]
for a grid function $V$ defined on the grid $\mathcal{T}_{\mathbf{h}}$.  
We then define the (second order) central discrete Laplacian operator by 
\begin{equation} \label{LaplacianH}
	\Delta_{\mathbf{h}} \equiv \sum_{i=1}^d \delta_{x_i, h_i}^2 . 
\end{equation}
In the formulation we will also consider the ``staggered" (second order) central difference operators 
$\delta_{x_i, 2h_i}^2$ and $\Delta_{2\mathbf{h}}$ defined by replacing $\mathbf{h}$ with $2\mathbf{h}$ in 
\eqref{fd_xx} and \eqref{LaplacianH}, respectively.  
The ``staggered" operators have a 5-point stencil in each Cartesian direction and are defined using 
3 nodes.  
Note that ``ghost-values"  need to be 
introduced in order for the staggered difference operators to be well-defined at interior nodes adjacent to 
the boundary of $\Omega$, 
i.e., points in $\left( \cT_{\bh} \cap \Omega \right) \setminus \mathring{\cT}_{\bh}$.


\subsection{The Lax-Friedrich's method and monotonicity}\label{LF_background_sec}

In this section we introduce the convergence framework of Crandall and Lions and the notion of monotonicity.  
We then introduce the Lax-Friedrich's method as an example of a 
method that falls within the Crandall and Lions framework.  
We will again consider the Lax-Friedrich's method in Section~\ref{numerics_sec} 
as a baseline for gauging the performance of our proposed method 
in a series of numerical tests.  

The (convergence) framework of Crandall and Lions relies upon two fundamental concepts:  
consistency and monotonicity.  
For ease of presentation, we will absorb the term $\theta u$ into the operator $H$ in \eqref{HJ_pde} 
and assume the FD scheme has the form 
\begin{equation}\label{FD_mon_scheme}
	\widehat{H} \left( \nabla_\mathbf{h}^+ U_\alpha , \nabla_\mathbf{h}^- U_\alpha , U_\alpha , x_\alpha \right) = 0 
	\qquad \text{for all } \bx_\alpha \in \cT_{\bh} \cap \Omega , 
\end{equation}
where $\widehat{H}$ is called a numerical Hamiltonian.  
Observe that $\widehat{H}$ depends upon two discrete gradients.  
Assuming $\widehat{H}$ is continuous, we say $\widehat{H}$ is {\it consistent} if 
$\widehat{H} \left( \mathbf{q} , \mathbf{q} , u , \bx \right) = H(\mathbf{q}, u , \bx)$ 
for any vector $\mathbf{q} \in \mathbb{R}^d$, $u \in \mathbb{R}$, and $\bx \in \Omega$.  
For the stationary problem, 
we say $\widehat{H}$ is {\it monotone} on $[-R,R]$ if it is 
nondecreasing with respect to the node $U_\alpha$ and 
nonincreasing with respect to each 
node $U_{\alpha^\prime}$ such that $\bx_{\alpha^\prime} \neq \bx_\alpha$ for $\bx_{\alpha'}$ in the local stencil 
centered at $\bx_\alpha$ 
whenever $\left| \delta_{x_i, h_i}^+ U_{\alpha^\prime} \right| \leq R$.  
Consequently, the scheme is monotone if it is nonincreasing with respect to $\nabla_{\mathbf{h}}^+ U_\alpha$ 
and nondecreasing with respect to $\nabla_{\mathbf{h}}^- U_\alpha$ and $U_\alpha$.  

A standard example of a consistent, monotone scheme is the Lax-Friedrich's method.  
The Lax-Friedrich's numerical Hamiltonian is defined by 
\begin{equation} \label{LF_scheme}
	\widehat{H}_{\text{LF}} [ U_\alpha ] 
	\equiv H \left( \frac12 \nabla_{\mathbf{h}}^+ U_\alpha + \frac12 \nabla_{\mathbf{h}}^- U_\alpha , 
			U_\alpha , \bx_\alpha \right) 
		- \vec{\beta} \cdot \left( \nabla_{\mathbf{h}}^+ - \nabla_{\mathbf{h}}^- \right) U_\alpha , 
\end{equation}
where $\vec{\beta} \geq \vec{0}$ using the natural partial ordering for vectors.  
Observe that, by the Lipschitz condition for $H$ and the assumption that $H$ is increasing with 
respect to $u$, 
if each component of $\vec{\beta}$ is sufficiently large, 
then the Lax-Friedrich's numerical Hamiltonian is nonincreasing with respect to 
$\nabla_{\mathbf{h}}^+ U_\alpha$ 
and nondecreasing with respect to $\nabla_{\mathbf{h}}^- U_\alpha$ and $U_\alpha$. 
Thus, the method is monotone for $\vec{\beta}$ sufficiently large.  
Under the convergence framework of Barles-Souganidis, 
the Lax-Friedrich's method will converge to the viscosity solution of \eqref{HJ} or, for $\epsilon > 0$ fixed 
and $-\epsilon \Delta_{\bh} U_\alpha$ added to the formulation, 
the viscosity solution of \eqref{HJ_vis}.  

The term $- \vec{\beta} \cdot \left( \nabla_{\mathbf{h}}^+ - \nabla_{\mathbf{h}}^- \right) U_\alpha$ is called a 
{\it numerical viscosity} due to the fact 
\begin{align*}
	- \vec{\beta} \cdot \left( \nabla_{\mathbf{h}}^+ - \nabla_{\mathbf{h}}^- \right) U_\alpha 
	& = - \sum_{i=1}^d \beta_i \left( \delta_{x_i, h_i}^+ - \delta_{x_i, h_i}^- \right) U_\alpha \\ 
	& = - \sum_{i=1}^d \beta_i \frac{U_{\alpha + \mathbf{e}_i} - 2 U_\alpha + U_{\alpha - \mathbf{e}_i}}{h_i} \\ 
	& = - \sum_{i=1}^d \beta_i h_i \delta_{x_i, h_i}^2 U_\alpha . 
\end{align*}
Letting $\vec{\beta} = \beta \vec{1}$ for some constant $\beta \geq 0$, 
we have $- \vec{\beta} \cdot \left( \nabla_{\mathbf{h}}^+ - \nabla_{\mathbf{h}}^- \right) U_\alpha$ 
is a (second-order) 
central difference approximation of $-\beta \sum_{i=1}^d u_{x_i x_i}(\bx_\alpha) = -\beta \Delta u(\bx_\alpha)$ 
scaled by $h$ when the mesh is uniform. 
Hence, the Lax-Friedrich's scheme for approximating \eqref{HJ} can be rewritten as 
\begin{subequations} \label{FD_LF}
\begin{alignat}{2}
	\widehat{H}_{\text{LF}} [ U_\alpha ] \equiv - \beta \sum_{i=1}^d h_i \delta_{x_i, h_i}^2 U_\alpha 
		+ H \left( \nabla_{\mathbf{h}} U_\alpha , U_\alpha , \bx_\alpha \right) + \theta U_\alpha & = 0 
		\qquad 
		&& \text{if } \bx_\alpha \in \mathcal{T}_{\mathbf{h}} \cap \Omega , \\ 
	U_\alpha - g(\bx_\alpha) & = 0 \quad && \text{if } \bx_\alpha \in \cT_{\bh} \cap \partial \Omega , 
\end{alignat}
\end{subequations} 
and we can see that the scheme is limited to first order accuracy due to the numerical viscosity.  
The scheme essentially chooses $\epsilon = \beta h$ in \eqref{HJ_vis} making it a direct 
realization of the vanishing viscosity approach.  

The first order bound is a consequence of the monotonicity.  
Every convergent monotone FD scheme for \eqref{HJ} implicitly approximates the differential equation 
\begin{equation} \label{tadmor_eqn}
	- \beta h ``\Delta u" + H(\nabla u , u ,\bx) = 0 
\end{equation}
for sufficiently large and possibly nonlinear $\beta \geq 0$, 
where $- \beta h ``\Delta u"$ corresponds to a numerical viscosity (c.f. \cite{Tadmor}).  
Thus, we cannot increase the power of the coefficient $h$ for the numerical viscosity 
and maintain a convergent, monotone scheme.  
To achieve a higher order scheme we will weaken the monotonicity condition and 
instead rely upon a higher-order stabilization technique 
that introduces a vanishing moment instead of a vanishing viscosity.


\section{Formulation and consistency analysis} \label{formulation_sec}

We approximate the solution of \eqref{HJ} using a FD method with up to 
second order local truncation errors.    
We will introduce two different choices for an auxiliary boundary condition that ensure the 
corresponding system of equations is well-defined.   
The main idea when formulating the new method 
will be to add a high-order stabilization term in lieu of the first-order accurate numerical viscosity 
in \eqref{tadmor_eqn}.  
The new stabilization term will sacrifice the monotonicity assumption
creating technicalities for the admissibility and stability analysis.  
We will formulate both a new method and a modified version that uses a cutoff operator 
as part of a corrector to the diffusive error in the Lax-Friedrich's method.  
The new method without a cutoff will be capable of second-order accuracy, 
and, for certain choices of the cutoff, will have the same solution as the modified version 
for which we can prove various admissibility, stability, and convergence results in Section~\ref{analysis_sec}.  
The new methods can trivially be extended for approximating solutions to \eqref{HJ_vis}.  


\subsection{Formulation and consistency analysis of a high-order corrector scheme} \label{unfiltered_sec}

The new high-order, non-monotone method is defined as follows. 
Let $\gamma > 0$, $\beta \geq 0$, and $p \in [0,1]$ be constants.  
Define the set of nodes $\mathcal{S}_{h_i}$ by 
\begin{align}\label{Sh_grid}
	\mathcal{S}_{h_i} \equiv \big\{ \bx_\alpha \in \cT_{\bh} \cap \partial \Omega \mid & \,  
		\bx_{\alpha} + h_i \mathbf{e}_i \in \cT_{\bh} \cap \Omega \text{ or } 
		\bx_{\alpha} - h_i \mathbf{e}_i \in \cT_{\bh} \cap \Omega \big\} 
\end{align}
for all $i \in \{1,2,\ldots,d\}$.  
The proposed high-order FD method for approximating solutions to \eqref{HJ} 
is defined as finding a grid function $U_\alpha$ such that 
\begin{subequations} \label{FD_scheme}
\begin{alignat}{2} 
	\widehat{H}_{\mathbf{h}} [ U_\alpha ]   & = 0 
		&& \text{if } \bx_\alpha \in \cT_{\bh} \cap \Omega , \label{FD_scheme1} \\ 
	U_\alpha - g(\bx_\alpha) & =0 \quad && \text{if } \bx_\alpha \in \cT_{\bh} \cap \partial \Omega , \label{FD_scheme2} \\ 
	B_{\bh} U_\alpha & = 0 \quad && \text{if } \mathbf{x}_\alpha\in \mathcal{S}_{h_i} \subset \mathcal{T}_{\mathbf{h}} \cap \partial\Omega  \label{FD_scheme3} 
\end{alignat} 
\end{subequations} 
for 
\begin{align*}
	\widehat{H}_{\mathbf{h}} [ U_\alpha ]  
	& \equiv - \beta h^2 \Delta_{\bh} U_\alpha 
		+ H \left( \nabla_{\mathbf{h}} U_\alpha , U_\alpha , \bx_\alpha \right) + \theta U_\alpha 
		+ \gamma h^p \left( \Delta_{2\mathbf{h}} - \Delta_{\mathbf{h}} \right) U_\alpha \\ 
	& \equiv H_{\bh} [U_\alpha] + \gamma h^p \left( \Delta_{2\mathbf{h}} - \Delta_{\mathbf{h}} \right) U_\alpha
\end{align*}
and the auxiliary boundary condition operator $B_{\bh}$ defined by either 
\begin{align} \label{bc2a}
	B_{\bh} U_\alpha \equiv - \delta_{x_i, h_i}^2 U_\alpha 
\end{align} 
or 
\begin{align} \label{bc2b}
	B_{\bh} U_\alpha \equiv  - \delta_{x_i, h_i}^2 U_\alpha + \delta_{x_i, h_i}^2 U_{\alpha'}  , 
\end{align} 
where 
$\bx_{\alpha'} \in \{ \bx_{\alpha} - h_i \mathbf{e}_i , \bx_{\alpha} + h_i \mathbf{e}_i \}$ 
such that $\bx_{\alpha'} \in \cT_{\bh} \cap \Omega$.  
Note that $\mathcal{S}_{h_i} \neq \emptyset$ for at most 
one $i \in \{1,2,\ldots,d\}$ at a given node $\bx_\alpha \in \cT_{\bh} \cap \partial \Omega$.  
The modified version defined in Section~\ref{filtered_sec} will use a cutoff operator that modifies the $\Delta_{2\bh}$ 
operator to ensure admissibility, stability, and convergence properties hold 
despite the non-monotone structure.   

In the formulation, we could use the auxiliary boundary condition 
$-\Delta_{\bh} U_\alpha = 0$ on $\cT_{\bh} \cap \partial \Omega$ 
instead of \eqref{bc2a} to be consistent with the choice in \cite{FDhjb}; 
however, \eqref{bc2a} is more closely matched with the chosen high-order 
stabilization term associated with $\gamma > 0$.
The newly proposed auxiliary boundary condition \eqref{bc2b} will allow higher-order accuracy.    
In general, we choose $\gamma$ sufficiently large so that the discrete operator 
$H_{\bh} - \frac{\gamma}{2} h^p \Delta_{\bh}$ is monotone.  
This is a mild assumption for $p < 1$ and comparable with choosing the numerical viscosity coefficient 
in the Lax-Friedrich's method for $p=1$.  
The choice $p=0$ is consistent with the Lax-Friedrich's-like method proposed and analyzed in 
\cite{FDhjb} for approximating uniformly elliptic fully nonlinear problems 
while the choice $p=1$ exploits the fact that the gradient terms for the Hamilton-Jacobi operator 
correspond to a $1/h$ scaling in order to achieve smaller truncation errors.  
We also note that, for $\beta > 0$, the scheme allows for a higher-order numerical viscosity to preserve 
the higher-order local truncation errors.  

\begin{remark}
Rewriting $\widehat{H}_{\bh}[U_\alpha] = - \gamma h^p \Delta_{\bh} U_\alpha + H_{\bh}[U_\alpha] 
+ \gamma h^p \Delta_{2 \bh} U_\alpha$, we can interpret the proposed method 
as adding a corrector term to the Lax-Friedrich's method when $\gamma$ is sufficiently large with $p=1$.  
The Lax-Friedrich's method is highly diffusive so the staggered discrete Laplacian offsets part of the 
numerical diffusivity.  
We can interpret the addition of the staggered discrete Laplacian operator as a corrector to account 
for the inconsistency with the Hamiltonian $H$ caused by adding the viscous term.  
The additional term will allow for higher order accuracy 
while the cutoff introduced in Section~\ref{filtered_sec} 
will guarantee admissibility, stability in $\ell^\infty$, and convergence 
for a modified corrector.   
\end{remark}

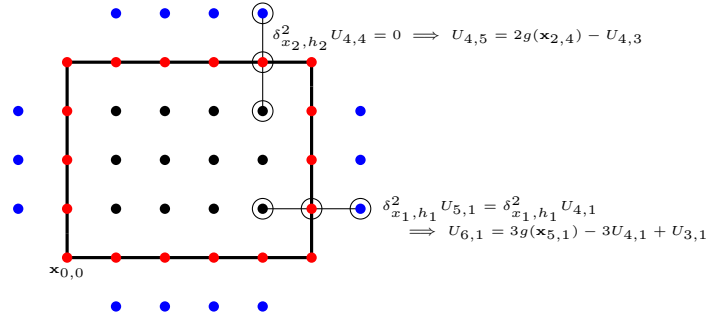
\begin{figure}
\begin{center}
\begin{tikzpicture}[scale=0.65] 
\draw[very thick] (0,-3.5) -- (0,0) -- (5,0) -- (5,-0.5); 
\draw[very thick] (0,-3.5) -- (0,-4) -- (5,-4) -- (5,-0.5);
\node[below] at (0,-4) {\tiny $\bx_{0,0}$};
\draw[thin] (4,-1) -- (4,1);
\node[right] at (4,0.5) {\tiny $\delta_{x_2, h_2}^2 U_{4,4} = 0 \implies U_{4,5} = 2 g(\bx_{2,4}) - U_{4,3}$};
\draw[thin] (4,-3) -- (6,-3);
\node[right] at (6.25,-3) {\tiny $\delta_{x_1, h_1}^2 U_{5,1} = \delta_{x_1, h_1}^2 U_{4,1}$}; 
	\node at (10,-3.5) {\tiny $\implies U_{6,1} = 3 g(\bx_{5,1}) - 3 U_{4,1} + U_{3,1}$};
\draw[thin] (4,1) circle (6pt);
\draw[thin] (4,0) circle (6pt);
\draw[thin] (4,-1) circle (6pt);
\draw[thin] (5,-3) circle (6pt);
\draw[thin] (4,-3) circle (6pt);
\draw[thin] (6,-3) circle (6pt);
\fill[color=red] (0,0) circle (3pt);
\fill[color=red] (1,0) circle (3pt);
\fill[color=red] (2,0) circle (3pt);
\fill[color=red] (3,0) circle (3pt);
\fill[color=red] (4,0) circle (3pt);
\fill[color=red] (5,0) circle (3pt);
\fill[color=red] (0,-1) circle (3pt);
\fill (1,-1) circle (3pt);
\fill (2,-1) circle (3pt);
\fill (3,-1) circle (3pt);
\fill (4,-1) circle (3pt);
\fill[color=red] (5,-1) circle (3pt);
\fill[color=red] (0,-2) circle (3pt);
\fill (1,-2) circle (3pt);
\fill (2,-2) circle (3pt);
\fill (3,-2) circle (3pt);
\fill (4,-2) circle (3pt);
\fill[color=red] (5,-2) circle (3pt);
\fill[color=red] (0,-3) circle (3pt);
\fill (1,-3) circle (3pt);
\fill (2,-3) circle (3pt);
\fill (3,-3) circle (3pt);
\fill (4,-3) circle (3pt);
\fill[color=red] (5,-3) circle (3pt);
\fill[color=red] (0,-4) circle (3pt);
\fill[color=red] (1,-4) circle (3pt);
\fill[color=red] (2,-4) circle (3pt);
\fill[color=red] (3,-4) circle (3pt);
\fill[color=red] (4,-4) circle (3pt);
\fill[color=red] (5,-4) circle (3pt);
\fill[color=blue] (1,1) circle (3pt);
\fill[color=blue] (2,1) circle (3pt);
\fill[color=blue] (3,1) circle (3pt);
\fill[color=blue] (4,1) circle (3pt);
\fill[color=blue] (-1,-1) circle (3pt);
\fill[color=blue] (-1,-2) circle (3pt);
\fill[color=blue] (-1,-3) circle (3pt);
\fill[color=blue] (6,-1) circle (3pt);
\fill[color=blue] (6,-2) circle (3pt);
\fill[color=blue] (6,-3) circle (3pt);
\fill[color=blue] (1,-5) circle (3pt);
\fill[color=blue] (2,-5) circle (3pt);
\fill[color=blue] (3,-5) circle (3pt);
\fill[color=blue] (4,-5) circle (3pt);
\end{tikzpicture}
\end{center}
\label{bc2_fig}
\caption{
A two-dimensional example of the mesh $\cT_{\bh}$ and enforcement of the auxiliary boundary condition \eqref{bc2a} 
or \eqref{bc2b}.  
The solid line corresponds to $\partial \Omega$.  
The black nodes are the unknown values in $\mathcal{T}_{\mathbf{h}} \cap \Omega$.  
The red nodes are the Dirichlet boundary data in $\mathcal{T}_{\mathbf{h}} \cap \partial \Omega$.  
The blue nodes represent ghost points and are uniquely determined by the choice of the auxiliary boundary condition.   
}
\end{figure}

The auxiliary boundary condition \eqref{FD_scheme3} is used to account for the 
fact that ghost points are needed to define $\widehat{H}_{\mathbf{h}} [ U_\alpha ]$ 
for all $\bx_\alpha \in \cT_{\bh} \cap \Omega$ one layer from $\partial \Omega$ due to the presence of the 
operator $\Delta_{2 \bh}$.  
Formally, the boundary condition extends the solution to the PDE along 
vectors normal to the boundary of $\Omega$ using either a linear or quadratic extension.  
See Figure~\ref{bc2_fig} for an example of the mesh, ghost points, and boundary conditions in two dimensions.  


We refer to the term $\gamma h^p \left( \Delta_{2\mathbf{h}} - \Delta_{\mathbf{h}} \right) U_\alpha$ 
as a numerical moment.  
A simple calculation reveals 
\begin{align} \label{moment_expansion}
	\left( \Delta_{2\mathbf{h}} - \Delta_{\mathbf{h}} \right) U_\alpha 
	& = \sum_{i=1}^d \left( \delta_{x_i, 2h_i}^2 - \delta_{x_i, h_i}^2 \right) U_\alpha 
	= \frac14 \sum_{i=1}^d h_i^2 \delta_{x_i, h_i}^2 \delta_{x_i, h_i}^2 U_\alpha \\ 
	\nonumber & = \frac14 \sum_{i=1}^d \left( \delta_{x_i, h_i}^2 U_{\alpha - \mathbf{e}_i} 
		+ \delta_{x_i, h_i}^2 U_{\alpha + \mathbf{e}_i} 
		- 2 \delta_{x_i, h_i}^2 U_\alpha \right) \\ 
	\nonumber & = \frac14 \sum_{i=1}^d h_i^2 \frac{U_{\alpha+2\mathbf{e}_i} - 4 U_{\alpha+\mathbf{e}_i} + 6 U_\alpha 
		- 4 U_{\alpha - \mathbf{e}_i} + U_{\alpha - 2 \mathbf{e}_i}}{h_i^4} . 
\end{align}
Then, for a smooth function $u \in C^6(\overline{\Omega})$, we have 
\[
	\left( \Delta_{2\mathbf{h}} - \Delta_{\mathbf{h}} \right) u(\bx_\alpha)
	= \frac14 \sum_{i=1}^d h_i^2 \left( u_{x_i x_i x_i x_i}(\bx_\alpha) + \frac16 h_i^2 \frac{\partial^6 u}{\partial x_i^6}(\xi_i^{(\alpha)}) \right) 
\]
for some $\xi_i^{(\alpha)} \in [\bx_\alpha - 2 h_i \mathbf{e}_i , \bx_\alpha + 2 h_i \mathbf{e}_i]$.  
The difference operator $\left( \Delta_{2\mathbf{h}} - \Delta_{\mathbf{h}} \right) U_\alpha$ is a (second-order) 
central difference approximation of the fourth order differential operator 
$\sum_{i=1}^d u_{x_i x_i x_i x_i}(\bx_\alpha)$ scaled by a 
constant proportional to $h^2$ when the mesh is quasi-uniform.  
As such, the proposed scheme is a direct realization of the vanishing moment method of Feng and Neilan.  

Suppose $\bx_\alpha \in \mathcal{S}_{h_i}$.  
A simple calculation reveals that the boundary condition \eqref{FD_scheme3} with 
$B_{\bh}$ defined by \eqref{bc2b} is equivalent to 
requiring 
\begin{align} \label{bc2b_ghost}
	& U_{\alpha+2 \mathbf{e}_i} - 3 U_{\alpha + \mathbf{e}_i} + 3 U_{\alpha} - U_{\alpha - \mathbf{e}_i} = 0 \\ 
	& \nonumber \qquad \text{or}  \\ 
	& \nonumber  
	U_{\alpha+\mathbf{e}_i} - 3 U_\alpha + 3 U_{\alpha - \mathbf{e}_i} - U_{\alpha - 2 \mathbf{e}_i} = 0 
\end{align}
depending on if $\bx_{\alpha + \mathbf{e}_i} \in \Omega$ or $\bx_{\alpha - \mathbf{e}_i} \in \Omega$.  
Thus, the auxiliary boundary condition is equivalent to assuming $u_{x_i x_i x_i} = 0$ 
using a (forward or backwards) first-order approximation of the third derivative enforced at the ghost point.  
By assuming the second derivative is constant through the boundary instead of zero-valued, 
the auxiliary boundary condition \eqref{bc2b} introduces a higher-order error for the boundary layer 
when compared to \eqref{bc2a}.
Formally, this is consistent with the fact that the choice \eqref{bc2a} would be exact if $u$ is linear while 
the choice \eqref{bc2b} 
would be exact if $u$ is quadratic.

\begin{lemma}\label{local_truncation_lemma}
The scheme \eqref{FD_scheme} has second order local truncation error over $\mathring{\cT_{\bh}}$.  
The choice of the auxiliary boundary operator \eqref{bc2a} gives a local truncation error of $\mathcal{O}(h^p)$ 
over $( \cT_{\bh} \cap \Omega ) \backslash\mathring{\cT_{\bh}}$, 
and the choice of the auxiliary boundary operator \eqref{bc2b} gives a local truncation error of  $\mathcal{O}(h^{p+1})$ over 
$( \cT_{\bh} \cap \Omega ) \backslash\mathring{\cT_{\bh}}$.  
\end{lemma}

\begin{proof}
Suppose $u(\bx) \in C^5(\overline{\Omega})$, and choose $\bx_\alpha \in \Omega$.  
Observe that, by the Lipschitz continuity of $H$, the mean value theorem, and the exact enforcement of the 
Dirichlet boundary condition, 
there holds 
\begin{align*}
	- \beta h^2 \Delta_{\bh} u(\bx_\alpha) 
		& = - \beta h^2 \sum_{i=1}^d \left( u_{x_i x_i}(\bx_\alpha) + \mathcal{O}(h_i^2) \right) \\ 
		& = \mathcal{O}(h^2) , \\   
	H \left( \nabla_{\bh} u(\bx_\alpha) , u(\bx_\alpha) , \bx_\alpha \right) 
		& = H \left( \nabla u(\bx_\alpha) + \mathcal{O}(h^2) \vec{1} , u(\bx_\alpha) , \bx_\alpha \right) \\ 
		& = \mathcal{O}(h^2) + H \left( \nabla u(\bx_\alpha) , u(\bx_\alpha) , \bx_\alpha \right) , \\ 
	\gamma h^p \left( \Delta_{2\mathbf{h}} - \Delta_{\mathbf{h}} \right) u (\bx_\alpha)  
		& = \gamma h^p \frac14 \sum_{i=1}^d h_i^2 \delta_{x_i, h_i}^2 \delta_{x_i, h_i}^2 u(\bx_\alpha) \\ 
		& = \frac{\gamma}4 h^p \sum_{i=1}^d h_i^2 \left( u_{x_i x_i x_i x_i}(\bx_\alpha) + \mathcal{O}(h) \right) \\ 
		& = \mathcal{O}(h^{2+p}) 
\end{align*}
for all $\bx_\alpha \in \mathring{\cT_{\bh}}$, 
where we have used the facts that the first and second bounds require $u \in C^4(\overline{\Omega})$ 
and the remainder term for the numerical moment requires $u \in C^5(\overline{\Omega})$.  
Plugging the above into \eqref{FD_scheme} for $0 \leq p \leq 1$, we have  
\[
	\widehat{H}_{\mathbf{h}} [u(\bx_\alpha)] = H \left( \nabla u(\bx_\alpha) , u(\bx_\alpha) , \bx_\alpha \right) 
	+ \theta U_\alpha 
	+ \mathcal{O}(h^2) 
\]
for all $\bx_\alpha \in \mathring{\cT_{\bh}}$.  

Suppose $\bx_\alpha \in (\cT_{\bh} \cap \Omega) \setminus \mathring{\cT_{\bh}}$, 
and note that the auxiliary boundary condition only has an impact on the local truncation error 
associated with the numerical moment while the other terms in \eqref{FD_scheme} would still have a second 
order local truncation error.  
Assume $\bx_{\alpha - \mathbf{e}_i} \in \partial \Omega$.  
Then, the choice \eqref{bc2a} implies 
\begin{align*}
   & \left( \Delta_{2\mathbf{h}} - \Delta_{\mathbf{h}} \right) u (\bx_\alpha)  \\ 
		& \qquad = - \frac14 \sum_{i=1}^d \delta_{x_i, h_i}^2 u(\bx_\alpha) \\ 
			& \qquad \qquad - \sum_{i=1}^d \frac{h_i}{4} \frac{ u(\bx_{\alpha} + 2 h_i \mathbf{e}_i) 
			- 3 u(\bx_{\alpha} + h_i \mathbf{e}_i) + 3 u(\bx_\alpha) - u(\bx_{\alpha} - h_i \mathbf{e}_i)}{h_i^3} \\
		& \qquad = - \frac14 \sum_{i=1}^d \left( u_{x_i x_i}(\bx_\alpha) 
			+ h_i u_{x_i x_i x_i}(\bx_{\alpha-\mathbf{e}_i}) + \mathcal{O}(h_i^2) \right) = \mathcal{O}(1) . 
\end{align*}
Since the choice \eqref{bc2b} implies 
$U_{\alpha - 2\mathbf{e}_i} -2U_{\alpha - \mathbf{e}_i}+U_{\alpha} 
= U_{\alpha - \mathbf{e}_i} -2U_{\alpha}+U_{\alpha + \mathbf{e}_i}$, 
there holds $3U_{\alpha } = 3 U_{\alpha - \mathbf{e}_i} - U_{\alpha - 2\mathbf{e}_i} +U_{\alpha + \mathbf{e}_i}$. 
Thus, for \eqref{bc2b}, there holds 
\begin{align*}
   & \left( \Delta_{2\mathbf{h}} - \Delta_{\mathbf{h}} \right) u (\bx_\alpha)  \\ 
		& \qquad = \sum_{i=1}^d \frac{h_i}{2} \frac{-u(\bx_{\alpha}-2h_i\mathbf{e}_i) +2 u(\bx_{\alpha}-h_i\mathbf{e}_i) 
			- 2 u(\bx_{\alpha} + h_i \mathbf{e}_i) + u(\bx_{\alpha} + 2 h_i \mathbf{e}_i)}{2h_i^3} \\
		& \qquad = \sum_{i=1}^d \frac{h_i}{2} \left( u_{x_i x_i x_i}(\bx_\alpha) + \mathcal{O}(h_i^2) \right) 
		= \mathcal{O}(h).
\end{align*}
The case when $\bx_{\alpha + \mathbf{e}_i} \in \partial \Omega$ is analogous.  
\hfill
\end{proof}

\begin{corollary}
The numerical moment is exact for linear functions when using the auxiliary boundary operator \eqref{bc2a}, 
and it is exact for quadratic functions when using the auxiliary boundary operator \eqref{bc2b}.  
For $\beta > 0$, the numerical viscosity scaled by $h^2$ presents an $\mathcal{O}(h^2)$ consistency error, 
and it is only exact for linear functions.  
The scheme is exact for linear functions for either choice of auxiliary boundary condition, 
and it is exact for quadratic functions when $\beta = 0$ and \eqref{bc2b} is used.  
\end{corollary}

We lastly look more closely at properties of $-\Delta_{2 \bh}$ and reformulate the numerical moment by removing the ghost values from the formulation.  
We show that the matrix representation of the minus staggered discrete Laplacian operator $-\Delta_{2\bh}$ 
is a nonsingular M-matrix 
when incorporating the auxiliary condition \eqref{FD_scheme3} with either \eqref{bc2a} or \eqref{bc2b} 
defining the auxiliary boundary operator.  
Assume $\bx_{\alpha - \mathbf{e}_i} \in \partial \Omega$.  
Then, the equation for $U_{\alpha}$ becomes 
\begin{align*}
	2 U_\alpha - U_{\alpha + 2 \mathbf{e}_i} - U_{\alpha - 2 \mathbf{e}_i} 
	& = 2 U_\alpha - U_{\alpha + 2 \mathbf{e}_i} - 2 U_{\alpha - \mathbf{e}_i} 
		+ U_{\alpha}  \\ 
	& = 3 U_\alpha - U_{\alpha + 2 \mathbf{e}_i} - 2 g(\bx_{\alpha-\mathbf{e}_i})
\end{align*}
when enforcing \eqref{FD_scheme3} with \eqref{bc2a} and, by \eqref{bc2b_ghost}, there holds 
\begin{align*}
	2 U_\alpha - U_{\alpha + 2 \mathbf{e}_i} - U_{\alpha - 2 \mathbf{e}_i} 
	& = 2 U_\alpha - U_{\alpha + 2 \mathbf{e}_i} 
		- 3 U_{\alpha - \mathbf{e}_i} + 3 U_{\alpha} - U_{\alpha + \mathbf{e}_i}  \\ 
	& = 5 U_\alpha - U_{\alpha + \mathbf{e}_i} - U_{\alpha + 2 \mathbf{e}_i} 
		- 3 g(\bx_{\alpha-\mathbf{e}_i})
\end{align*}
when enforcing \eqref{FD_scheme3} with \eqref{bc2b}.  
We have similar results if $\bx_{\alpha + \mathbf{e}_i} \in \partial \Omega$, 
and it follows that the matrix representation of $-\Delta_{2 \bh}$ is a nonsingular M-matrix.  
Plugging these values into the numerical moment, we have 
\begin{align*}
	\delta_{x_i, 2h_i}^2 U_\alpha - \delta_{x_i, h_i}^2 U_\alpha 
	& = \frac{-3 U_\alpha + U_{\alpha + 2 \mathbf{e}_i} + 2 g(\bx_{\alpha-\mathbf{e}_i})}{4 h_i^2} 
		+ \frac{2 U_\alpha - g(\bx_{\alpha-\mathbf{e}_i}) - U_{\alpha + \mathbf{e}_i}}{h_i^2} \\ 
	& = \frac{5 U_\alpha - 4 U_{\alpha + \mathbf{e}_i} + U_{\alpha + 2 \mathbf{e}_i} 
		- 2 g(\bx_{\alpha-\mathbf{e}_i}) }{4 h_i^2}
\end{align*}
when enforcing \eqref{FD_scheme3} with \eqref{bc2a} and 
\begin{align*}
	\delta_{x_i, 2h_i}^2 U_\alpha - \delta_{x_i, h_i}^2 U_\alpha 
	& = \frac{-5 U_\alpha + U_{\alpha + \mathbf{e}_i} + U_{\alpha + 2 \mathbf{e}_i} 
		+ 3 g(\bx_{\alpha-\mathbf{e}_i})}{4 h_i^2} \\ 
		& \qquad 
		+ \frac{2 U_\alpha - g(\bx_{\alpha-\mathbf{e}_i}) - U_{\alpha + \mathbf{e}_i}}{h_i^2} \\ 
	& = \frac{3 U_\alpha - 3 U_{\alpha + \mathbf{e}_i} + U_{\alpha + 2 \mathbf{e}_i} 
		- g(\bx_{\alpha-\mathbf{e}_i}) }{4 h_i^2}
\end{align*}
when enforcing \eqref{FD_scheme3} with \eqref{bc2b} for all $i=1,2,\ldots,d$.
Thus, the numerical moment is still increasing with respect to the nodes at $\bx_\alpha$ and 
$\bx_{\alpha \pm 2 \mathbf{e}_i}$ and decreasing with respect to the nodes at $\bx_{\alpha \pm \mathbf{e}_i}$ 
when removing the ghost points from the formulation.  
This type of monotonicity is similar to the generalized monotonicity in \cite{Kellie} 
and prevents the application of standard analytic techniques for monotone methods.


\subsection{Formulation of the modified high-order corrector scheme} \label{filtered_sec}

In this section, we formulate a modified version of the FD method \eqref{FD_scheme} that 
ensures the non-monotone corrector associated with the numerical moment 
does not lead to an instability in the $\ell^\infty$ norm.  
The method will require choosing upper and lower cutoff functions.  
We will discuss how to choose the cutoff functions in Section~\ref{admissiblity_stability_sec} 
to guarantee convergence to the viscosity solution of \eqref{HJ}.  
In the numerical experiments in Section~\ref{numerics_sec}, the modified scheme often 
yields a solution to \eqref{FD_scheme} with second-order accuracy.  

We first motivate the choice for the cutoff operators.  
Suppose $\bx_\alpha \in \cT_{\bh} \cap \Omega$. 
By the definition of the scheme \eqref{FD_scheme} and rewriting the numerical moment, we have 
\begin{align*}
	\widehat{H}_{\bh} [ V_\alpha ] 
	& = - \beta h^2 \Delta_{\bh} V_\alpha + H(\nabla_{\bh} V_\alpha , V_\alpha , \bx_\alpha) + \theta V_\alpha 
		+ \gamma h^p \left( \Delta_{2 \bh} - \Delta_{\bh} \right) V_\alpha \\ 
	& = - \beta h^2 \Delta_{\bh} V_\alpha + H(\nabla_{\bh} V_\alpha , V_\alpha , \bx_\alpha) + \theta V_\alpha \\ 
		& \qquad 
		+ \gamma h^p \sum_{i=1}^d \frac{ V_{\alpha+2\mathbf{e}_i} - 4 V_{\alpha+\mathbf{e}_i} + 6 V_\alpha 
		- 4 V_{\alpha - \mathbf{e}_i} + V_{\alpha - 2 \mathbf{e}_i }}{4 h_i^2} \\ 
	& = - \beta h^2 \Delta_{\bh} V_\alpha + H(\nabla_{\bh} V_\alpha , V_\alpha , \bx_\alpha) + \theta V_\alpha \\ 
		& \qquad 
		+ \gamma h^p \sum_{i=1}^d \frac{ - 2 V_{\alpha+\mathbf{e}_i} + 6 V_\alpha 
			- 2 V_{\alpha - \mathbf{e}_i} }{4 h_i^2} \\ 
		& \qquad 
		+ \gamma h^p \sum_{i=1}^d 
			\frac{ V_{\alpha+2\mathbf{e}_i} - 2 V_{\alpha+\mathbf{e}_i} 
				+ V_{\alpha-2\mathbf{e}_i} - 2 V_{\alpha-\mathbf{e}_i} }{4 h_i^2} \\ 
	& \equiv H_{\bh}[V_\alpha] 
		- \gamma \frac12 h^p \Delta_{\bh} V_\alpha + \gamma h^p \frac12 \sum_{i=1}^d \frac{1}{h_i^2} V_\alpha
		- \gamma h^p \frac14 \sum_{i=1}^d \frac{1}{h_i^2} \left( L_{i}^+ V_\alpha + L_{i}^- V_\alpha \right) 
\end{align*}
for 
\[
	L_{i}^\pm V_\alpha \equiv  2 V_{\alpha \pm \mathbf{e}_i} - V_{\alpha \pm 2 \mathbf{e}_i} 
\]
for all $i=1,2,\ldots,d$.  
The cutoffs will be based on placing upper and lower bounds on the difference operators $L_i^\pm$.  

Choose functions $\underline{U}, \overline{U} \in S(\cT_{\bh})$ such that $\underline{U} \leq \overline{U}$, 
and define the truncation operators $\overline{\underline{L}}_i^\pm$ by 
\[
	\overline{\underline{L}}_i^\pm V_\alpha \equiv \begin{cases}
	\overline{U}_\alpha \qquad & \text{if } L_i^\pm V_\alpha > \overline{U}_\alpha , \\ 
	\underline{U}_\alpha \qquad & \text{if } L_i^\pm V_\alpha < \underline{U}_\alpha , \\ 
	L_i^\pm V_\alpha \qquad & \text{otherwise} .  
	\end{cases}
\]
Then, the modified version of the proposed FD method \eqref{FD_scheme} for approximating solutions to \eqref{HJ} 
is defined as finding a grid function $U_\alpha$ such that 
\begin{subequations} \label{FD_filtered}
\begin{alignat}{2} 
	\widehat{\overline{\underline{H}}}_{\mathbf{h}} [ U_\alpha ]   & = 0 
		&& \text{if } \bx_\alpha \in \cT_{\bh} \cap \Omega , \label{FD_filtered1} \\ 
	U_\alpha - g(\bx_\alpha) & =0 \quad && \text{if } \bx_\alpha \in \cT_{\bh} \cap \partial \Omega , \label{FD_filtered2} \\ 
	B_{\bh} U_\alpha & = 0 \quad && \text{if } \mathbf{x}_\alpha\in \mathcal{S}_{h_i} \subset \mathcal{T}_{\mathbf{h}} \cap \partial\Omega  \label{FD_filtered3} 
\end{alignat} 
\end{subequations} 
for 
\begin{align*}
	\widehat{\overline{\underline{H}}}_{\mathbf{h}} [ U_\alpha ]  
	& \equiv H_{\bh}[U_\alpha] 
		- \gamma \frac12 h^p \Delta_{\bh} U_\alpha + \gamma h^p \frac12 \sum_{i=1}^d \frac{1}{h_i^2} U_\alpha
		- \gamma h^p \frac14 \sum_{i=1}^d \frac{1}{h_i^2} \left( 
			\overline{\underline{L}}_{i}^+ U_\alpha + \overline{\underline{L}}_{i}^- U_\alpha \right) \\ 
	& = \widehat{\overline{\underline{H}}}_{\text{LF}}[U_\alpha]
		+ \gamma h^p \frac12 \sum_{i=1}^d \frac{1}{h_i^2} U_\alpha
		- \gamma h^p \frac14 \sum_{i=1}^d \frac{1}{h_i^2} \left( 
			\overline{\underline{L}}_{i}^+ U_\alpha + \overline{\underline{L}}_{i}^- U_\alpha \right) 
\end{align*}
with $\widehat{\overline{\underline{H}}}_{\text{LF}} \equiv H_{\bh} - \frac{\gamma}{2} h^p \Delta_{\bh}$. 
Note that $\widehat{\overline{\underline{H}}}_{\text{LF}}$ is a monotone Lax-Friedrich's operator 
for $\gamma$ sufficiently large even when $\beta = 0$ in the definition of $H_{\bh}$.  
If $U$ defined by \eqref{FD_filtered} satisfies $\overline{\underline{L}}_i^\pm U_\alpha = L_i^\pm U_\alpha$ 
for all $\bx_\alpha \in \cT_{\bh} \cap \Omega$, then \eqref{FD_filtered} reduces to \eqref{FD_scheme} 
and has the potential for higher-order consistency errors.  
In practice, we choose $\overline{U} = U_{\text{LF}} + C h$ and $\underline{U} = U_{\text{LF}} - Ch$ for 
$U_{\text{LF}}$ the solution to $\widehat{\overline{\underline{H}}}_{\text{LF}} [U_\alpha] = 0$ with 
$U_\alpha = g(\bx_\alpha)$ for all $\bx_\alpha \in \cT_{\bh} \cap \partial \Omega$ 
and $C > 0$.  
See Section~\ref{admissiblity_stability_sec} for more details 
and Section~\ref{numerics_sec} for numerical tests corresponding to this choice.


\section{Analytic Results} \label{analysis_sec}

In this section, we introduce a new analytic technique for proving admissibility, stability, and convergence 
results for the modified high-order corrector scheme \eqref{FD_filtered} in the non-degenerate case when $\theta > 0$.
The admissibility analysis will naturally extend to the degenerate case $\theta = 0$ by adding a 
numerical stabilization term $h^2 U_\alpha$, but the technique for proving the uniform stability 
bounds will no longer directly apply.  

The admissibility and stability analysis will be based on choosing appropriate discrete sub- and supersolutions that 
uniformly bound the FD approximation.  
Let $S(\cT_{\bh})$ denote the space of all (grid) functions mapping $\cT_{\bh}$ to $\mathbb{R}$.  
The analysis will apply the Schauder fixed-point theorem to show that the proposed FD scheme \eqref{FD_filtered} has a solution.  
To this end, observe that the space $S(\cT_{\bh})$ paired with the $\ell^\infty$ metric is a Banach space.  
Let $\underline{U}, \overline{U} \in S(\cT_{\bh})$ such that $\underline{U} \leq \overline{U}$, 
and define the set $\overline{\underline{S}}(\cT_{\bh})$ by 
\[
	\overline{\underline{S}}(\cT_{\bh}) \equiv \left\{ V \in S(\cT_{\bh}) \mid \underline{U}_\alpha \leq V_\alpha \leq 
	\overline{U}_\alpha \text{ for all } \bx_\alpha \in \cT_{\bh} \right\} . 
\]
Then, $\overline{\underline{S}}(\cT_{\bh})$ is a nonempty, convex, compact subset of $S(\cT_{\bh})$, 
and the Schauder fixed-point theorem ensures that every continuous function that maps 
$\overline{\underline{S}}(\cT_{\bh})$ into itself has a fixed point.  
We will choose $\underline{U}$ and $\overline{U}$ appropriately 
and formulate a continuous mapping whose fixed point is a solution to the proposed FD scheme.  
The stability result will follow by the choice of $\underline{U}$ and $\overline{U}$, 
where $\underline{U}$ will be an appropriate subsolution and $\overline{U}$ will be an appropriate supersolution.  

In this section, we will first apply the general framework that will be used to analyze the modified 
high-order corrector scheme to the Lax-Friedrich's method in Section~\ref{LF_analysis-sec} 
to establish the techniques and some preliminary results.  
We next motivate the use of cutoffs while exploring the impact of the non-monotone numerical moment stabilizer 
in Section~\ref{motivation_sec} before applying the ideas to the proposed scheme \eqref{FD_filtered} in Section~\ref{admissiblity_stability_sec}.  
Some final observations regarding the analysis framework and 
the proposed high-order corrector scheme \eqref{FD_scheme} can be found in Section~\ref{admissiblity_stability_unfiltered_sec}.

\subsection{Analysis of the Lax-Friedrich's scheme} \label{LF_analysis-sec}

We first consider established techniques and results for monotone FD methods \eqref{FD_mon_scheme} 
for approximating \eqref{HJ} 
to motivate the new techniques that we will employ. 
To this end, we will reference the Lax-Friedrich's method \eqref{FD_LF} 
and typical fixed-point reformulations of the discrete problem.  
Many of the ideas in this section are linked to the techniques used for proving existence  
results for nonlinear reaction-diffusion equations as discussed in \cite{LMZ}. 
We will derive admissibility, stability, uniqueness, and convergence rates for the Lax-Friedrich's method 
using Schauder fixed-point theory. 

Choose $\tau > 0$, 
and consider the mapping $\cM_{\tau} : S(\mathcal{T}_{\mathbf{h}}) \to S(\cT_{\bh})$ 
defined by $\widehat{V} \equiv \cM_\tau V$ if 
\begin{subequations} \label{Mtau}
\begin{alignat}{2} 
	& \widehat{V}_\alpha = V_\alpha - \tau \widehat{H}_{\text{LF}} [ V_\alpha ] \quad 
		&& \text{if } \bx_\alpha \in \mathcal{T}_{\mathbf{h}} \cap \Omega , \\ 
	& \widehat{V}_\alpha = g(\bx_\alpha) \quad && \text{if } \bx_\alpha \in \cT_{\bh} \cap \partial \Omega 
\end{alignat}
\end{subequations} 
for all $V \in S(\cT_{\bh})$.  
Then, for all $\tau$ sufficiently small, the mapping $\cM_\tau$ is monotone in the sense that 
the right hand side $V_\alpha - \tau \widehat{H}_{\text{LF}} [ V_\alpha ]$ is nondecreasing with respect to 
each component of $V$ when $H$ is globally Lipschitz 
(with ways to handle the case when $H$ is only locally Lipschitz).  
Furthermore, the mapping is a contraction in $\ell^\infty$ for $\theta > 0$.  
The fact that $\cM_\tau$ is a contraction mapping for $\tau$ sufficiently small 
guarantees the admissibility of the Lax-Friedrich's scheme as well as uniform stability estimates
based inversely on $\theta$ and 
directly on $\| H(\mathbf{0}, 0, \cdot) \|_{C^0(\Omega)}$ and $\| g \|_{C^0(\partial \Omega)}$.  
We derive similar results using a weaker fixed-point analysis.  

The admissibility and stability analysis will be based on choosing appropriate discrete sub- and supersolutions that 
uniformly bound the Lax-Friedrich's approximation.  
We say $\underline{U} \in S(\cT_{\bh})$ is a subsolution of \eqref{FD_LF} if 
\begin{alignat*}{2} 
	\widehat{H}_{\text{LF}} [ \underline{U}_\alpha ]   & \leq 0 
		&& \text{if } \bx_\alpha \in \cT_{\bh} \cap \Omega ,  \\ 
	\underline{U}_\alpha - g(\bx_\alpha) & \leq 0 \quad && \text{if } \bx_\alpha \in \cT_{\bh} \cap \partial \Omega .
\end{alignat*} 
We say $\overline{U} \in S(\cT_{\bh})$ is a supersolution of \eqref{FD_LF} if 
the inequalities are all changed to `$\geq$'.  
The goal will be to show \eqref{FD_LF} has a solution $U$ with 
$\underline{U} \leq U \leq \overline{U}$ for a suitable choice of 
subsolution $\underline{U}$ and supersolution $\overline{U}$. 

Observe that $\widehat{H}_{\text{LF}}[ c 1_{\bh}(\bx_\alpha)] = H(\mathbf{0},c,\bx_\alpha) + \theta c$ 
is nonnegative for $c$ sufficiently large and nonpositive for $c$ sufficiently small for $\theta > 0$, 
where $1_{\bh} \in S(\cT_{\bh})$ denotes the grid function such that $1_{\bh}(\bx_\alpha) = 1$ 
for all $\bx_\alpha \in \cT_{\bh}$.  
Assuming $|c| \geq \| g \|_{C^0(\partial \Omega)}$, 
we have there exists a value $c > 0$ independent of $h$ such that 
$\overline{U} = c 1_{\bh}$ is a discrete supersolution of \eqref{FD_LF}
and $\underline{U} = - c 1_{\bh}$ is a discrete subsolution of \eqref{FD_LF} 
since the numerical viscosity is zero-valued when applied to a constant-valued grid function.    
By the monotonicity of $\cM_\tau$ for all $\tau > 0$ sufficiently small 
when restricted to inputs in $\overline{\underline{S}}(\cT_{\bh})$ and using the fact that $H$ is locally Lipschitz
as well as the facts that $\underline{U}$ is a subsolution and 
$\overline{U}$ is a supersolution, 
if $\underline{U} \leq V \leq \overline{U}$, 
then $\widehat{V}_\alpha = V_\alpha - \tau \widehat{H}_{\text{LF}}[V_\alpha]$ satisfies   
\begin{align*}
	\underline{U}_\alpha 
	\leq \underline{U}_\alpha - \tau \widehat{H}_{\text{LF}}[ \underline{U}_\alpha] 
	\leq V_\alpha - \tau \widehat{H}_{\text{LF}}[V_\alpha] 
	\leq \overline{U}_\alpha - \tau \widehat{H}_{\text{LF}}[ \overline{U}_\alpha ] 
	\leq \overline{U}_\alpha 
\end{align*}
for all $\bx_\alpha \in \cT_{\bh} \cap \Omega$ and 
$\underline{U}_\alpha \leq \widehat{V}_\alpha = g(\bx_\alpha) \leq \overline{U}_\alpha$ 
for all $\bx_\alpha \in \cT_{\bh} \cap \partial \Omega$.  
Hence, $\underline{U} \leq \widehat{V} \leq \overline{U}$, and it follows that 
\eqref{FD_LF} has a solution $U$ with $\underline{U} \leq U \leq \overline{U}$ by the 
Schauder fixed-point theorem.  
Furthermore, the bound $c$ for the fixed-point $U$ scales inversely with respect to $\theta$ and 
directly with respect to $\| H(\mathbf{0}, 0, \cdot) \|_{C^0(\Omega)}$ and $\| g \|_{C^0(\partial \Omega)}$.  

The  admissibility and stability framework for the Lax-Friedrich's method \eqref{FD_LF} 
can be used to derive a uniqueness result.  
We will then use the consistency of the scheme to derive a convergence result.  
\begin{theorem}
The solution to \eqref{FD_LF} is unique 
when choosing $\beta > 0$ sufficiently large to ensure 
$\widehat{H}_{\text{LF}}$ is monotone 
for any grid function in $\overline{\underline{S}}(\cT_{\bh})$.
\end{theorem}

\begin{proof}
First, suppose the PDE operator $H$ is linear.  
Then, the Lax-Friedrich's method corresponding to approximating 
the linear problem 
\begin{subequations} \label{HJ_linear}
\begin{alignat}{2}
	\vec{b} \cdot \nabla u + c u + \theta u & = f && \text{in } \Omega , \\ 
	u & = g \qquad && \text{on } \partial \Omega 
\end{alignat}
\end{subequations}
with $f, \vec{b}, c$ bounded over $\Omega$ and $c \geq 0$
has a solution.  
Since any linear problem corresponding to the FD method 
has a solution for all possible data corresponding to $f$ 
and $g$, the corresponding matrix must be nonsingular.  
Thus, the FD solution is unique when $H$ is linear.  

Suppose the FD method \eqref{FD_LF} 
applied to \eqref{HJ} with $H$ nonlinear has two solutions $U, V \in \overline{\underline{S}}(\cT_{\bh})$.  
Let $W = U - V$.  
Then, by the Lipchitz continuity of $H$ and the mean value theorem, there exists a linear operator $L_{\bh}$ such that 
$L_{\bh} [W_\alpha] = \widehat{H}_{\text{LF}}[U_\alpha] - \widehat{H}_{\text{LF}}[V_\alpha]$.  
In particular, $L_{\bh}[W_\alpha]$ would correspond to the linear operator 
\begin{align}\label{linearization}
	L_{\bh} [W_\alpha] 
	& = - \beta h \Delta_{\bh} W_\alpha 
		+ \vec{b}(\bx_\alpha) \cdot \nabla_{\bh} W_\alpha + c (\bx_\alpha) W_\alpha
		+ \theta W_\alpha 
\end{align}
for some bounded functions $\vec{b}, c$ with $c \geq 0$.  
Furthermore, there holds $W_\alpha = 0$ for all $\bx_\alpha \in \cT_{\bh} \cap \partial \Omega$ 
and $L_{\bh} [W_\alpha] = 0$ for all $\bx_\alpha \in \cT_{\bh} \cap \Omega$.  
Since a matrix representation of $L_{\bh}$ must be nonsingular, it follows that $W_\alpha = 0$ 
for all $\bx_\alpha \in \cT_{\bh}$.  
\hfill 
\end{proof}

\begin{theorem} \label{LF_convergence_thm}
Let $u$ be the solution to \eqref{HJ}, and assume $u$ is sufficiently smooth so that the local truncation error 
for the FD scheme \eqref{FD_LF} is $\mathcal{O}(h^r)$, 
where $r \leq 1$ and $r=1$ if $u \in C^2(\overline{\Omega})$.  
Let $U$ be the solution to \eqref{FD_LF} 
with $\beta \geq 0$ sufficiently large to ensure $\widehat{H}_{\text{LF}}$
is monotone for any grid function in $\overline{\underline{S}}(\cT_{\bh})$.
Then, $\| u - U \|_{\ell^\infty(\cT_{\bh})} = \mathcal{O}(h^r)$.  
\end{theorem}

\begin{proof}
There exists a uniformly bounded grid function $C$ such that 
$\widehat{H}_{\text{LF}}[u(\bx_\alpha)] - \widehat{H}_{\text{LF}}[U_\alpha] = C_\alpha h^r$ 
for all $\bx_\alpha \in \cT_{\bh} \cap \Omega$ with $u(\bx_\alpha) - U_\alpha = 0$ 
for all $\bx_\alpha \in \cT_{\bh} \cap \partial \Omega$.  
Let $E_\alpha = u (\bx_\alpha) - U_\alpha$ for all $\bx_\alpha \in \cT_{\bh}$.  
Then, since $H$ is locally Lipschitz, by the mean value theorem there exists a linear operator $L_{\bh}$ such that 
$L_{\bh}[E_\alpha] = C_\alpha h^r$, where $L_{\bh}$ has the form given by \eqref{linearization}.  
Using the constant-valued sub- and supersolutions 
$\overline{U} = \frac{2 \| C \|_{\ell^\infty(\cT_{\bh})}}{\theta} h^r 1_{\bh}$ 
and 
$\underline{U} = -\frac{2 \| C \|_{\ell^\infty(\cT_{\bh})}}{\theta} h^r 1_{\bh}$, 
we have  
$\| E \|_{\ell^\infty(\cT_{\bh})} \leq \frac{2 \| C \|_{\ell^\infty(\cT_{\bh})}}{\theta} h^r$ 
applying the Schauder fixed-point theorem, 
and the result follows. 
\hfill 
\end{proof}

\begin{remark}
The convergence result can be extended to lower-regularity viscosity solutions by the Barles-Souganidis framework 
for admissible, stable, consistent, and monotone methods which guarantees locally uniform convergence 
of the Lax-Friedrich's approximation.  
\end{remark}

\subsection{Analytic motivation for the modified high-order formulation} \label{motivation_sec}

The mapping $\cM_\tau$ defined by \eqref{Mtau} is not monotone for any choice $\tau > 0$ when replacing $\widehat{H}_{LF}$ 
with $\widehat{H}_{\bh}$ from the proposed high-order scheme in Section~\ref{unfiltered_sec} since the mapping would be 
decreasing with respect to the $U_{\alpha \pm 2 \mathbf{e}_i}$ nodes introduced by the corrector term in the numerical moment.  
Consequently, we cannot use the local update to directly bound $\widehat{V}_\alpha$ by $\| V \|_{\ell^\infty(\cT_{\bh})}$.  
To see this, suppose $V_\alpha = \overline{U}_\alpha$, 
$V_{\alpha \pm \mathbf{e}_i} = \overline{U}_{\alpha \pm \mathbf{e}_i}$, and 
$V_{\alpha \pm 2 \mathbf{e}_i} = \underline{U}_{\alpha \pm 2 \mathbf{e}_i}$ 
for all $i=1,2,\ldots,d$ 
for some $\bx_\alpha \in \mathring{\cT}_{\bh} \cap \Omega$, 
where the sub- and supersolution are the constant-valued functions chosen in Section~\ref{LF_analysis-sec}.  
Then, by \eqref{moment_expansion} and using the fact the sub- and supersolution are constant-valued with 
$\underline{U} = - \overline{U}$, there holds 
\begin{align*}
	\widehat{V}_\alpha & = V_\alpha - \tau \widehat{H}_{\bh} [ V_\alpha ] \\ 
	& = \left( 1 - \tau \theta - \tau \sum_{i=1}^d \frac{\beta 4 h^2 + \gamma 3 h^p}{2 h_i^2} \right) \overline{U}_\alpha 
		- \tau H (\mathbf{0} , \overline{U}_\alpha , \bx_\alpha) \\ 
		& \qquad 
		+ 2 \tau \sum_{i=1}^d \left( \beta \frac{h^2}{h_i^2} + \gamma \frac{h^p}{h_i^2} \right) 
			\overline{U}_{\alpha} 
		- \tau \gamma \sum_{i=1}^d \frac{h^p}{2 h_i^2} \underline{U}_\alpha \\ 
	& = \overline{U}_\alpha - \tau H (\mathbf{0} , \overline{U}_\alpha , \bx_\alpha) 
		- \tau \left( \theta - \sum_{i=1}^d \frac{\gamma h^p}{h_i^2} \right) \overline{U}_\alpha .  
\end{align*}
Therefore, $\widehat{V}_\alpha > \overline{U}_\alpha$ independent of $\tau$ 
if $h$ is small enough such that  
\[
	\gamma \sum_{i=1}^d \frac{h^p}{h_i^2} \overline{U}_\alpha 
	> H(\mathbf{0}, \overline{U}_\alpha , \bx_\alpha) + \theta \overline{U}_\alpha \geq 0 . 
\]
Since $p \in [0,1]$, we would expect $\widehat{V}_\alpha > \overline{U}_\alpha$ for reasonably fine meshes when 
$\gamma > 0$ is fixed.  
This type of issue naturally arises from the inclusion of the numerical moment 
since it does not have sign control based on a monotone structure like the numerical viscosity.  
The numerical moment allows extreme positive or negative terms 
based on the distance between the sub- and supersolutions which cannot directly be controlled by the 
sub- and supersolutions themselves.  

An approach we can use is to modify the troublesome terms in the 
mapping so that they are {\em eventually} constant-valued (and hence monotone).  
Thus, we can use truncated difference operators that control the size of the numerical moment.  
Choose $M > 0$, and let $\psi(x) = - \sum_{i=1}^d \frac{M}{2d} \left( x - \frac{a_i+b_i}{2} \right)^2 + C$ 
with $C$ chosen such that 
$\psi$ restricted to $\cT_{\bh} \cap \overline{\Omega}$ is a supersolution of 
\eqref{FD_scheme} 
and $-\psi$ restricted to $\cT_{\bh} \cap \overline{\Omega}$ is a subsolution of 
\eqref{FD_scheme}.  
Note that $-\Delta_{2 \bh} \psi = - \Delta_{\bh} \psi = M$.  
Define the truncated difference operator $\widetilde{\Delta}_{2 \bh}$ by 
\[
\widetilde{\Delta}_{2 \bh} V_\alpha \equiv 
\begin{cases}
	-\Delta_{2 \bh} \psi(\bx_\alpha) & \text{ if }  \Delta_{2 \bh} V_\alpha > M , \\ 
	\Delta_{2 \bh} V_\alpha & \text{ if } -M \leq \Delta_{2 \bh} V_\alpha \leq M , \\ 
	\Delta_{2 \bh} \psi(\bx_\alpha) & \text{ if } \Delta_{2 \bh} V_\alpha < - M 
\end{cases}
\]
and the modified numerical operator $\widetilde{H}_{\bh}$ by 
$\widetilde{H}_{\bh} [V_\alpha] 
\equiv H_{\bh} [V_\alpha] + \gamma h^p (\widetilde{\Delta}_{2 \bh} - \Delta_{\bh}) V_\alpha$
for all $V \in S(\cT_{\bh})$.    
Then, the operator is eventually constant-valued with respect to the $U_{\alpha \pm 2 \mathbf{e}_i}$ 
nodes that break the standard monotonicity condition.  

Consider the mapping $\widetilde{\cM}_{\tau} : S(\mathcal{T}_{\mathbf{h}}) \to S(\cT_{\bh})$ 
defined by $\widehat{V} \equiv \widetilde{\cM}_\tau V$ if 
\begin{alignat*}{2} 
	& \widehat{V}_\alpha = V_\alpha - \tau \widetilde{H}_{\bh} [ V_\alpha ] \quad 
		&& \text{if } \bx_\alpha \in \mathcal{T}_{\mathbf{h}} \cap \Omega , \\ 
	& \widehat{V}_\alpha - g(\bx_\alpha) = 0 \quad && \text{if } \bx_\alpha \in \cT_{\bh} \cap \partial \Omega , \\ 
	& B_{\bh} \widehat{V}_\alpha = B_{\bh} V_\alpha \quad && \text{if } \bx_\alpha \in \mathcal{S}_{h_i} 
\end{alignat*}
for all $V \in S(\cT_{\bh})$.  
Suppose $- \psi(\bx_\alpha) \leq V_\alpha \leq \psi(\bx_\alpha)$ for all $\bx_\alpha \in \cT_{\bh}$.  
Then, for $\gamma$ sufficiently large such that 
the difference operator 
$H_{\bh} - \gamma h^p \Delta_{\bh}$ is monotone 
and $\tau$ sufficiently small, 
there holds 
\begin{align*}
	\widehat{V}_\alpha 
	& = V_\alpha - \tau \gamma  h^p \widetilde{\Delta}_{2 \bh} V_\alpha 
		- \tau \left( H_{\bh}[V_\alpha] - \gamma h^p \Delta_{\bh} V_\alpha \right) \\ 
	& \leq V_\alpha - \tau \gamma  h^p \widetilde{\Delta}_{2 \bh} \psi(\bx_\alpha) 
		- \tau \left( H_{\bh}[V_\alpha] - \gamma h^p \Delta_{\bh} V_\alpha \right) \\ 
	& \leq \psi(\bx_\alpha) - \tau \gamma  h^p \widetilde{\Delta}_{2 \bh} \psi(\bx_\alpha) 
		- \tau \left( H_{\bh}[\psi(\bx_\alpha)] - \gamma h^p \Delta_{\bh} \psi(\bx_\alpha) \right) \\ 
	& = \psi(\bx_\alpha) - \tau \widehat{H}_{\bh} [ \psi(\bx_\alpha) ] \\ 
	& \leq \psi(\bx_\alpha) 
\end{align*}
for all $\bx_\alpha \in \cT_{\bh} \cap \Omega$.  
Similarly, $\widehat{V}_\alpha \geq - \psi(\bx_\alpha)$, 
and it follows that the truncated problem with $\widehat{H}_{\bh}$ replaced by $\widetilde{H}_{\bh}$ in \eqref{FD_scheme} has a solution $U$ with $-\psi \leq U \leq \psi$.  
Furthermore, the solution solves the unmodifed problem if there exists an $M$ sufficiently large such that 
$\widetilde{\Delta}_{2 \bh} U_\alpha = \Delta_{2 \bh} U_\alpha$ 
for all $\bx_\alpha \in \cT_{\bh} \cap \Omega$.  
In the next section we will instead use the modified high-order correction scheme \eqref{FD_filtered} 
as a way to truncate the numerical moment.  
The alternative cutoff function used in \eqref{FD_filtered} will allow for a strategic choice of the sub- and supersolutions 
that will guarantee uniform stability and convergence results.  

\begin{remark}
In the counterexample for directly replacing $\widehat{H}_{LF}$ with $\widehat{H}_{\bh}$ from the proposed scheme, 
the value for $\widehat{V}_\alpha$ can remain in the interval $[\underline{U}_\alpha , \overline{U}_\alpha]$ by making $\tau$ 
sufficiently small for a given input $V$ where $V_\alpha \in ( \underline{U}_\alpha , \overline{U}_\alpha )$.  
Necessarily $\tau \to 0^+$ as $V_\alpha$ approaches the bounds if the numerical moment has the wrong sign.  
\end{remark}

\begin{remark}
Sending $\tau \to 0^+$, the iteration $\widetilde{\cM}_\tau$ can be regarded as a pseudo-timestepping 
algorithm using ``concavity limiters" as opposed to slope limiters.  
As such, the ideas in this paper should have applications to deriving stability bounds for higher-order 
(non-monotone) methods for dynamic Hamilton-Jacobi equations.  
\end{remark}

\subsection{Analysis of the modified high-order corrector scheme} \label{admissiblity_stability_sec}

In this section we will first analyze the admissibility and stability of the FD scheme \eqref{FD_filtered}.  
We will then choose particular sub- and supersolutions that guarantee convergence.  
The analysis will heavily rely upon the results for the Lax-Friedrich's method in Section~\ref{LF_analysis-sec}.  
Notationally, we let $U_{\text{LF}}$ be the Lax-Friedrich's approximation corresponding to solving 
\begin{subequations} \label{FD_LFgamma}
\begin{alignat}{2} 
	\widehat{\overline{\underline{H}}}_{\text{LF}} [ U_\alpha ]   & = 0 
		&& \text{if } \bx_\alpha \in \cT_{\bh} \cap \Omega , \label{FD_LFgamma1} \\ 
	U_\alpha - g(\bx_\alpha) & =0 \quad && \text{if } \bx_\alpha \in \cT_{\bh} \cap \partial \Omega . \label{FD_LFgamma2}  
\end{alignat} 
\end{subequations} 
We also let $\overline{U}$ be a supersolution to \eqref{FD_LFgamma} and $\underline{U}$ be a subsolution to \eqref{FD_LFgamma} 
with $\underline{U} \leq \overline{U}$. 
Particular choices for $\overline{U}$ and $\underline{U}$ will be specified in the convergence analysis.  

The analysis is again based on the pseudo-timestepping iteration.  
Choose $\tau > 0$, 
and consider the mapping $\overline{\underline{\cM}}_{\tau} : S(\mathcal{T}_{\mathbf{h}}) \to S(\cT_{\bh})$ 
defined by $\widehat{V} \equiv \overline{\underline{\cM}}_\tau V$ if 
\begin{subequations} \label{Mtau_filtered}
\begin{alignat}{2} 
	& \widehat{V}_\alpha = V_\alpha - \tau \widehat{\overline{\underline{H}}}_{\bh} [ V_\alpha ] \quad 
		&& \text{if } \bx_\alpha \in \mathcal{T}_{\mathbf{h}} \cap \Omega , \\ 
	& \widehat{V}_\alpha = g(\bx_\alpha) \quad && \text{if } \bx_\alpha \in \cT_{\bh} \cap \partial \Omega , \\ 
	& B_{\bh} \widehat{V}_\alpha = B_{\bh} V_\alpha \quad && \text{if } \bx_\alpha \in \mathcal{S}_{h_i} 
\end{alignat}
\end{subequations} 
for all $V \in S(\cT_{\bh})$.  
Clearly a fixed-point of \eqref{Mtau_filtered} is a solution to \eqref{FD_filtered}.  
We show $\overline{\underline{\cM}}_{\tau}$ maps $\overline{\underline{S}}(\cT_{\bh})$ into itself.  

Suppose $\gamma > 0$ is sufficiently large to ensure $\widehat{\overline{\underline{H}}}_{\text{LF}}$ is monotone 
and $\tau > 0$ is small enough such that $V_\alpha - \tau \widehat{\overline{\underline{H}}}_{\bh} [ V_\alpha ]$ 
is nondecreasing with respect to $V_\alpha$ 
for all $\bx_\alpha \in \cT_{\bh} \cap \Omega$ for all $V \in \overline{\underline{S}}(\cT_{\bh})$.  
Let $V \in \overline{\underline{S}}(\cT_{\bh})$, and choose $\bx_\alpha \in \cT_{\bh} \cap \Omega$.  
Observe that, by the monotonicity of $\widehat{\overline{\underline{H}}}_{\text{LF}}$, 
by the definition of $\overline{\underline{L}}_{i}^\pm$, 
and since $\overline{U}$ is a supersolution of \eqref{FD_LFgamma}, 
there holds 
\begin{align*}
	& V_\alpha - \tau \widehat{\overline{\underline{H}}}_{\bh} [ V_\alpha ] \\ 
	& \qquad = V_\alpha 
		- \tau \widehat{\overline{\underline{H}}}_{\text{LF}}[V_\alpha]
		- \tau \gamma h^p \frac12 \sum_{i=1}^d \frac{1}{h_i^2} V_\alpha
		+ \tau \gamma h^p \frac14 \sum_{i=1}^d \frac{1}{h_i^2} \left( 
			\overline{\underline{L}}_{i}^+ V_\alpha + \overline{\underline{L}}_{i}^- V_\alpha \right) \\ 
	& \qquad \leq \overline{U}_\alpha 
		- \tau \widehat{\overline{\underline{H}}}_{\text{LF}}[\overline{U}_\alpha]
		- \tau \gamma h^p \frac12 \sum_{i=1}^d \frac{1}{h_i^2} \overline{U}_\alpha
		+ \tau \gamma h^p \frac14 \sum_{i=1}^d \frac{1}{h_i^2} \left( 
			\overline{\underline{L}}_{i}^+ V_\alpha + \overline{\underline{L}}_{i}^- V_\alpha \right) \\ 
	& \qquad \leq \overline{U}_\alpha 
		- \tau \widehat{\overline{\underline{H}}}_{\text{LF}}[\overline{U}_\alpha]
		- \tau \gamma h^p \frac12 \sum_{i=1}^d \frac{1}{h_i^2} \overline{U}_\alpha
		+ \tau \gamma h^p \frac12 \sum_{i=1}^d \frac{1}{h_i^2} \overline{U}_\alpha \\ 
	& \qquad = \overline{U}_\alpha 
		- \tau \widehat{\overline{\underline{H}}}_{\text{LF}}[\overline{U}_\alpha] \\ 
	& \qquad \leq \overline{U}_\alpha . 
\end{align*}
Similarly, there holds $V_\alpha - \tau \widehat{\overline{\underline{H}}}_{\bh} [ V_\alpha ] \geq \underline{U}_\alpha$.  
Thus, $\widehat{V} = \overline{\underline{\cM}}_\tau V \in \overline{\underline{S}}(\cT_{\bh})$, 
and it follows that \eqref{Mtau_filtered} has a fixed point in $\overline{\underline{S}}(\cT_{\bh})$ 
by the Schauder fixed-point theorem.  
We have proved the following.  

\begin{theorem} 
Let $\overline{U}$ be a supersolution to \eqref{FD_LFgamma} and $\underline{U}$ be a subsolution to \eqref{FD_LFgamma} 
with $\underline{U} \leq \overline{U}$ 
with $\gamma > 0$ sufficiently large to ensure $\widehat{\overline{\underline{H}}}_{\text{LF}}$ is monotone 
over $\overline{\underline{S}}(\cT_{\bh})$.  
Then \eqref{FD_filtered} has a solution in $\overline{\underline{S}}(\cT_{\bh})$.  
\end{theorem}

The choice of the cutoff functions determines the cutoff operators $\overline{\underline{L}}_{i}^\pm$ 
and how well the solution to \eqref{FD_filtered} approximates the viscosity solution to \eqref{HJ}.  
Note that $U_{\text{LF}}$ is both a sub- and supersolution to \eqref{FD_LFgamma}.  
Suppose the viscosity solution $u \in C^2(\overline{\Omega})$.  
Then, by Theorem~\ref{LF_convergence_thm} and the local truncation error for the Lax-Friedrich's scheme, 
we have there exists $C > 0$ depending on $D^2 u$ such that 
$\| u - U_{\text{LF}} \|_{\ell^\infty(\cT_{\bh})} \leq C h$.  
We also have $U_{\text{LF}} + \sigma 1_{\bh}$ is a supersolution of \eqref{FD_LFgamma}
and $U_{\text{LF}} - \sigma 1_{\bh}$ is a subsolution of \eqref{FD_LFgamma} 
for any constant $\sigma > 0$.  
Thus, we choose $\overline{U} = U_{\text{LF}} + c h 1_{\bh}$ and $\underline{U} = U_{\text{LF}} - c h 1_{\bh}$ 
for any constant $c > 0$ 
to ensure the solution to \eqref{FD_filtered} converges to the viscosity solution of \eqref{HJ} by the squeeze theorem.  
Note that, even if $u$ does not have sufficient regularity to ensure $\| u - U_{\text{LF}} \|_{\ell^\infty(\cT_{\bh})} \leq C h$, 
then $U_{\text{LF}}$ still converges to $u$ by the Barles-Souganidis theorem.  
Thus, the solution to \eqref{FD_filtered} will converge to the viscosity solution of \eqref{HJ} even when 
the solution has lower regularity.  
Choosing $c \gtrsim C$, we expect the solution to the modified high-order corrector scheme may be more accurate than the 
underlying Lax-Friedrich's approximation.  

\begin{theorem} \label{FD_filtered_convergence_thm}
Let $\overline{U} = U_{\text{LF}} + c h 1_{\bh}$ and $\underline{U} = U_{\text{LF}} - c h 1_{\bh}$ 
for a constant $c > 0$. 
Suppose $\gamma > 0$ is sufficiently large to ensure $\widehat{\overline{\underline{H}}}_{\text{LF}}$ is monotone 
over $\overline{\underline{S}}(\cT_{\bh})$. 
Then, the FD approximations solving \eqref{FD_filtered} converge to the viscosity solution of \eqref{HJ}.  
\end{theorem}

\begin{remark}
The modified high-order corrector scheme with the choice 
$\overline{U} = U_{\text{LF}} + c h 1_{\bh}$ and $\underline{U} = U_{\text{LF}} - c h 1_{\bh}$
acts like a predictor-corrector method that uses the Lax-Friedrich's method to find a coarse approximation and then 
seeks a corrected solution in a ball of radius $ch$.  
We expect the accuracy to increase as $c$ approaches or exceeds 
the constant $C$ based on the consistency error of the 
Lax-Friedrich's method.  
The modified scheme ensures the lack of monotonicity does not lead to an instability 
while allowing for potentially second order accuracy when the underlying viscosity solution is sufficiently smooth.  
The boost in accuracy is observed in all numerical tests in Section~\ref{numerics_sec}.  
\end{remark}

\begin{remark}
The analysis can easily be extended to non-uniform meshes. 
The choice of $\overline{U}$ and $\underline{U}$ can also be generalized to 
$\overline{U} = U_{\text{LF}} + c_{\bh}$ and $\underline{U} = U_{\text{LF}} - c_{\bh}$ for any sequence $c_{\bh} \to 0$ 
without impacting the various results.  
\end{remark}

\subsection{Applications to the high-order corrector scheme} \label{admissiblity_stability_unfiltered_sec}

The admissibility and stability analysis in Section~\ref{admissiblity_stability_sec} does not directly apply 
to the high-order method \eqref{FD_scheme} 
due to the lack of monotonicity for the underlying method 
and not using cutoffs to control the non-monotone contributions.  
However, the solution to \eqref{FD_filtered} may also solve \eqref{FD_scheme} 
which can easily be checked in post-processing.  
We see in all of the numerical tests that the solution to \eqref{FD_filtered} also solves \eqref{FD_scheme} 
when choosing the cutoffs in Theorem~\ref{FD_filtered_convergence_thm} for $c$ sufficiently large 
with the solution to \eqref{FD_scheme} having increased accuracy. 
Thus, the modified scheme \eqref{FD_filtered} appears to be capable of the higher-order accuracy linked with 
the proposed high-order scheme \eqref{FD_scheme} when choosing appropriate cutoffs while also 
having guaranteed admissibility, stability, and convergence properties.  
Future work seeks to directly analyze \eqref{FD_scheme} 
without having to use the modified scheme to account for the lack of monotonicity. 

\begin{figure}
\begin{center}
\begin{tikzpicture}[scale=0.65] 
\draw[very thick,<->] (-1,0) -- (8,0);
\node[right] at (8,0) {\tiny $x$};
\draw[very thick,<->] (0,-2) -- (0,7);
\node[above] at (0,7) {\tiny $y$};
\draw[thin] (2,-0.2) -- (2,0.2);
\node[below] at (2,-0.2) {\tiny $x_{\alpha-2}$};
\draw[thin] (3,-0.2) -- (3,0.2);
\node[below] at (3,-0.2) {\tiny $x_{\alpha-1}$};
\draw[thin] (4,-0.2) -- (4,0.2);
\node[below] at (4,-0.2) {\tiny $x_{\alpha}$};
\draw[thin] (5,-0.2) -- (5,0.2);
\node[below] at (5,-0.2) {\tiny $x_{\alpha+1}$};
\draw[thin] (6,-0.2) -- (6,0.2);
\node[below] at (6,-0.2) {\tiny $x_{\alpha+2}$};
\draw[thin] (-0.2,-1) -- (0.2,-1);
\node[left] at (-0.2,-1) {\tiny $\underline{U}_\alpha$};
\draw[thin] (-0.2,6) -- (0.2,6);
\node[left] at (-0.2,6) {\tiny $\overline{U}_\alpha$};
\fill (4,-1) circle (3pt);
\fill[color=red] (4,6) circle (5pt);
\fill[color=red] (2,3) circle (5pt);
\fill (2,3) circle (3pt);
\fill[color=red] (6,6) circle (5pt);
\fill (6,6) circle (3pt);
\fill (3,1) circle (3pt);
\fill[color=red] (3,4.5) circle (5pt);
\fill[color=red] (5,6) circle (5pt);
\fill (5,2.5) circle (3pt);
\draw[thin, color=red, densely dashed] (2,3) -- (4,6);
\draw[thin, color=red, densely dashed] (4,6) -- (6,6);
\draw[thin, color=black, densely dashed] (2,3) -- (4,-1);
\draw[thin, color=black, densely dashed] (4,-1) -- (6,6);
\draw[thin, color=blue, densely dashed] (3,1) -- (3,4.5);
\draw[thin, color=blue, densely dashed] (5,2.5) -- (5,6);
\draw[thin, color=cyan, densely dashed] (4,-1) -- (4,6);
\end{tikzpicture}
\end{center}
\label{local_V_pic}
\caption{
Sample plots locally bounding the values of $U_{\alpha \pm \mathbf{e}_i}$ to ensure $U$ solves 
both \eqref{FD_scheme} and \eqref{FD_filtered}.  
The values for $U_{\alpha \pm 2 \mathbf{e}_i}$ are fixed.  
The value for $U_\alpha$ can range over the entire interval $[\underline{U}_\alpha , \overline{U}_\alpha]$ 
which is shaded in cyan.  
The value for $U_{\alpha \pm \mathbf{e}_i}$ can range over the blue subintervals 
that are determined by the values of $U_{\alpha \pm 2\mathbf{e}_i}$, $\overline{U}_\alpha$, and $\underline{U}_\alpha$.  
}
\end{figure}
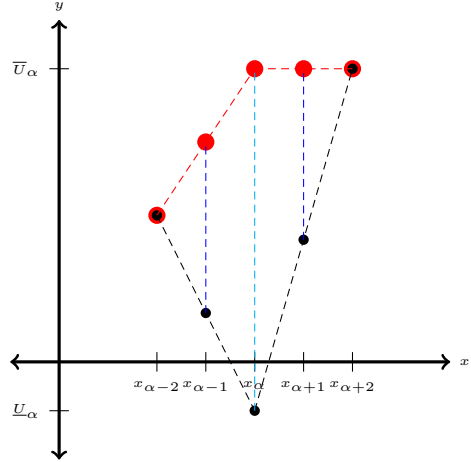

We now interpret the choice for the cutoffs when modifying the operators $L_i^\pm$ 
that were used in lieu of the truncated difference operators proposed in Section~\ref{motivation_sec}.  
Observe that $\underline{\overline{L}}_i^\pm U_\alpha = L_i^\pm U_\alpha$ if 
$\underline{U}_\alpha \leq L_i^\pm U_\alpha \leq \overline{U}_\alpha$ for all $\bx_\alpha \in \cT_{\bh} \cap \Omega$, 
where $U$ is a solution of \eqref{FD_filtered} with the boundary data and ghost values defined by 
\eqref{FD_filtered2} and \eqref{FD_filtered3}.  
Then, 
\[
	\underline{U}_\alpha \leq 2 U_{\alpha \pm \mathbf{e}_i} - U_{\alpha \pm 2 \mathbf{e}_i} \leq \overline{U}_\alpha 
	\implies 
	\frac{\underline{U}_\alpha + U_{\alpha \pm 2 \mathbf{e}_i}}{2} \leq U_{\alpha \pm \mathbf{e}_i} 
	\leq \frac{\overline{U}_\alpha + U_{\alpha \pm 2 \mathbf{e}_i}}{2} . 
\]
Thus, the possible values for $U_{\alpha \pm \mathbf{e}_i}$ are bounded away from extreme values 
based on the values for $U_{\alpha \pm 2 \mathbf{e}_i}$.  
The bounds allow $U_\alpha$ to be any value in $[\underline{U}_\alpha , \overline{U}_\alpha]$.  
If \eqref{FD_scheme} is admissible and stable, then 
$\underline{U}_\alpha + U_{\alpha \pm 2 \mathbf{e}_i}$ can be made arbitrarily small and 
$\overline{U}_\alpha + U_{\alpha \pm 2 \mathbf{e}_i}$ can be made arbitrarily large 
since $U_{\alpha \pm 2 \mathbf{e}_i}$ would be bounded, 
$\underline{U} = U_{\text{LF}} - c 1_{\bh}$ is an appropriate subsolution for any $c > 0$, 
and $\overline{U} = U_{\text{LF}} + c 1_{\bh}$ is an appropriate supersolution for any $c > 0$.  
See Figure~\ref{local_V_pic} for a visual representation of the bounds enforced by the modified operators 
$\underline{\overline{L}}_i^\pm U_\alpha$.  
We see that the modified operators restrict the concavity of $U$ when $U_\alpha$ is close to the upper bound 
$\overline{U}_\alpha$ or the lower bound $\underline{U}_\alpha$.  
This is in contrast to the motivational method $\widetilde{H}_{\bh}$ in Section~\ref{motivation_sec} 
that uniformly bounded the concavity instead of only locally bounding concavity for approximations near the 
upper and lower bounds. 

\section{Numerical experiments} \label{numerics_sec}

In this section we test the convergence rates of the proposed high-order corrector methods \eqref{FD_scheme} and \eqref{FD_filtered} 
for both choices of the auxiliary boundary condition \eqref{bc2a} and \eqref{bc2b} 
when approximating both smooth solutions and lower regularity viscosity solutions of \eqref{HJ}. 
We choose $\overline{U} = U_{\text{LF}} + c h 1_{\bh}$ and $\underline{U} = U_{\text{LF}} - c h 1_{\bh}$ 
when implementing \eqref{FD_filtered} for various choices of $c > 0$.  
For large values of $c$, we see that the solution to \eqref{FD_filtered} is also a solution to \eqref{FD_scheme}.  
We also benchmark the Lax-Friedrich's method \eqref{FD_LFgamma} and see that the proposed 
methods are always more accurate with \eqref{FD_scheme} capable of second-order accuracy.  
Errors will be measured in the $\ell^\infty$-norm. 
All meshes are uniform with $h_x = h_y$ for two-dimensional problems. 
In all tests we set $\beta = 0$ and $p=1$.  
The one-dimensional tests use the domain $(-1,1)$ and $\gamma = 10$, 
and the two-dimensional tests use the domain $(-1,1) \times (-1,1)$ and $\gamma = 5$.  

The tests are all performed in {\it Matlab} and use {\it fsolve} to solve the resulting nonlinear algebraic problem.  
The initial guess for the Lax-Friedrich's method is the zero function, and the 
initial guess for the proposed schemes is the Lax-Friedrich's approximation.  
When recording the results, the column labeled cutoff for the modified high-order corrector scheme indicates 
whether $\overline{\underline{L}}_i^\pm [U_\alpha] = L_i^\pm U_\alpha$ for all $\bx_\alpha \in \cT_{\bh} \cap \Omega$ 
or not.  
A value of `yes' indicates that the cutoff was applied.  
A value of `no' indicates that the modified scheme agreed with the proposed high-order corrector scheme.  
In all tests, a value of $c=10$ for the modified scheme \eqref{FD_filtered} yielded the solution to the proposed 
high order corrector scheme \eqref{FD_scheme}.  
In general, the accuracy of the modified method increased as $c$ increased.  
More tests for the scheme \eqref{FD_scheme} can be found in \cite{FD_CDR} where the impact of only satisfying 
the boundary condition in the viscosity sense is also explored.  


\subsection{One-dimensional tests} \label{1D_numerics_sec}

We first consider a series of numerical experiments in one dimension to test the accuracy of the proposed schemes. 
The examples are globally Lipschitz and satisfy the analytic assumptions made for the differential operator $H$.  
In the first example, we use a linear problem and track the behavior for various values of $c$.  
The linear problem allows for extremely fine meshes to better track the asymptotic behavior of the modified schemes with respect 
to the cutoff functions.  
The second example highlights the fact that the Lax-Friedrich's method is overly diffusive when compared to 
the proposed high-order correction schemes that use a numerical moment instead of a numerical viscosity.  
The proposed method \eqref{FD_scheme} always has rates of convergence at least as high as predicted 
by the consistency analysis.  


\subsubsection{Example 1:  linear operator with a smooth solution}
Consider the problem 
\[
	H[u] = (3x^2-x+4)u_x + (x^2 + 1)u - f(x) 
\]
with $f$ and the boundary data chosen such that the exact solution is $u(x) = x^3 + \cos(4x)$.   
The results for the Lax-Friedrich's method and the proposed high-order method can be found in 
Table~\ref{1D-1a}.  
The results for the modified high-order corrector method with various choices for $c$ can be found in Tables~\ref{1D-1b} and \ref{1D-1c}.  
The approximation corresponding to \eqref{FD_filtered} agrees with the approximation for \eqref{FD_scheme} for all 
tested values $c \geq 7$.

\begin{table}[htb] 
{\scriptsize 
\begin{center}
\begin{tabular}{| c | c | c | c | c |}
		\hline
	& \multicolumn{2}{|c|}{Lax-Friedrich's} & \multicolumn{2}{|c|}{Scheme \eqref{FD_scheme}} \\ 
		\hline
	 $h$ & Error & Order & Error & Order \\ 
		\hline
	2.02e-02 & 1.24e-01 &  & 2.13e-03 &  \\ 
		\hline
	6.69e-03 & 4.15e-02 & 0.99 & 2.20e-04 & 2.05 \\ 
		\hline
	3.34e-03 & 2.08e-02 & 1.00 & 5.40e-05 & 2.02 \\ 
		\hline
	2.00e-03 & 1.25e-02 & 1.00 & 1.93e-05 & 2.01 \\ 
		\hline
	1.00e-03 & 6.23e-03 & 1.00 & 4.80e-06 & 2.01 \\ 
		\hline
	5.00e-04 & 3.11e-03 & 1.00 & 1.20e-06 & 2.00 \\ 
		\hline
	4.00e-04 & 2.49e-03 & 1.00 & 7.65e-07 & 2.00 \\ 
		\hline
	2.00e-04 & 1.25e-03 & 1.00 & 1.91e-07 & 2.00 \\ 
		\hline
\end{tabular}
\end{center}
}
\caption{
Approximations for Example 1 in one dimension for the Lax-Friedrich's method and \eqref{FD_scheme} with boundary condition \eqref{bc2a}.  
Similar results hold for the boundary condition \eqref{bc2b}.
}
\label{1D-1a}
\end{table}

\begin{table}[htb] 
{\scriptsize 
\begin{center}
\begin{tabular}{| c | c | c | c | c | c | c | c | c | c |}
		\hline
	& \multicolumn{3}{|c|}{c = 1} & \multicolumn{3}{|c|}{c = 2} & \multicolumn{3}{|c|}{c = 4} \\ 
		\hline
	 $h$ & Cutoff & Error & Order & Cutoff & Error & Order & Cutoff & Error & Order \\ 
		\hline
	2.02e-02 & yes & 1.04e-01 & & yes & 8.40e-02 &  & yes & 4.39e-02 &  \\ 
		\hline
	6.69e-03 & yes & 3.48e-02 & 0.99 & yes & 2.81e-02 & 0.99 & yes & 1.48e-02 & 0.98 \\ 
		\hline
	3.34e-03 & yes & 1.74e-02 & 1.00 & yes & 1.41e-02 & 1.00 & yes & 7.41e-03 & 0.99 \\
		\hline
	2.00e-03 & yes & 1.05e-02 & 1.00 & yes & 8.45e-03 & 1.00 & yes & 4.45e-03 & 1.00 \\
		\hline
	1.00e-03 & yes & 5.23e-03 & 1.00 & yes & 4.23e-03 & 1.00 & yes & 2.22e-03 & 1.00 \\ 
		\hline
	5.00e-04 & yes & 2.61e-03 & 1.00 & yes & 2.11e-03 & 1.00 & yes & 1.11e-03 & 1.00 \\ 
		\hline
	4.00e-04 & yes & 2.09e-03 & 1.00 & yes & 1.69e-03 & 1.00 & yes & 8.90e-04 & 1.00 \\ 
		\hline
	2.00e-04 & yes & 1.05e-03 & 1.00 & yes & 8.45e-04 & 1.00 & yes & 4.45e-04 & 1.00 \\ 
		\hline
\end{tabular}
\end{center}
}
\caption{
Approximations for Example 1 in one dimension for \eqref{FD_filtered} with boundary condition \eqref{bc2a} 
using various smaller values for $c$ that ensure the cutoffs are applied.  
Similar results hold for the boundary condition \eqref{bc2b}.
}
\label{1D-1b}
\end{table}

\begin{table}[htb] 
{\scriptsize 
\begin{center}
\begin{tabular}{| c | c | c | c | c | c | c | c | c | c |}
		\hline
	& \multicolumn{3}{|c|}{c = 6.22} & \multicolumn{3}{|c|}{c = 6.226} & \multicolumn{3}{|c|}{c = 6.228} \\ 
		\hline
	 $h$ & Cutoff & Error & Order & Cutoff & Error & Order & Cutoff & Error & Order \\ 
		\hline
	2.02e-02 & yes & 2.02e-03 &  & yes & 2.02e-03 &  & yes & 2.02e-03 &  \\ 
		\hline
	6.69e-03 & yes & 2.11e-04 & 2.04 & yes & 2.11e-04 & 2.04 & yes & 2.11e-04 & 2.04 \\ 
		\hline
	3.34e-03 & yes & 5.22e-05 & 2.01 & yes & 5.22e-05 & 2.01 & yes & 5.27e-05 & 2.00 \\ 
		\hline
	2.00e-03 & yes & 1.89e-05 & 1.98 & yes & 1.87e-05 & 2.01 & no & 1.93e-05 & 1.97 \\ 
		\hline
	1.00e-03 & yes & 7.27e-06 & 1.38 & yes & 4.66e-06 & 2.00 & no & 4.80e-06 & 2.01 \\ 
		\hline
	5.00e-04 & yes & 3.08e-06 & 1.24 & yes & 1.19e-06 & 1.97 & no & 1.20e-06 & 2.00 \\ 
		\hline
	4.00e-04 & yes & 2.39e-06 & 1.14 & no & 7.65e-07 & 1.99 & no & 7.65e-07 & 2.00 \\ 
		\hline
	2.00e-04 & yes & 1.22e-06 & 0.98 & no & 1.91e-07 & 2.00 & no & 1.91e-07 & 2.00 \\ 
		\hline
\end{tabular}
\end{center}
}
\caption{
Approximations for Example 1 in one dimension for \eqref{FD_filtered} with boundary condition \eqref{bc2a} 
using various values for $c$ corresponding for the transition region where 
the modified scheme agrees with the proposed high-order scheme.  
Similar results hold for the boundary condition \eqref{bc2b}.
}
\label{1D-1c}
\end{table}


\subsubsection{Example 2:  nonlinear Lipschitz operator with a non-smooth solution} \label{test3_sec}

\begin{figure}[htb] 
\begin{center}
\includegraphics[width=0.29\textwidth]{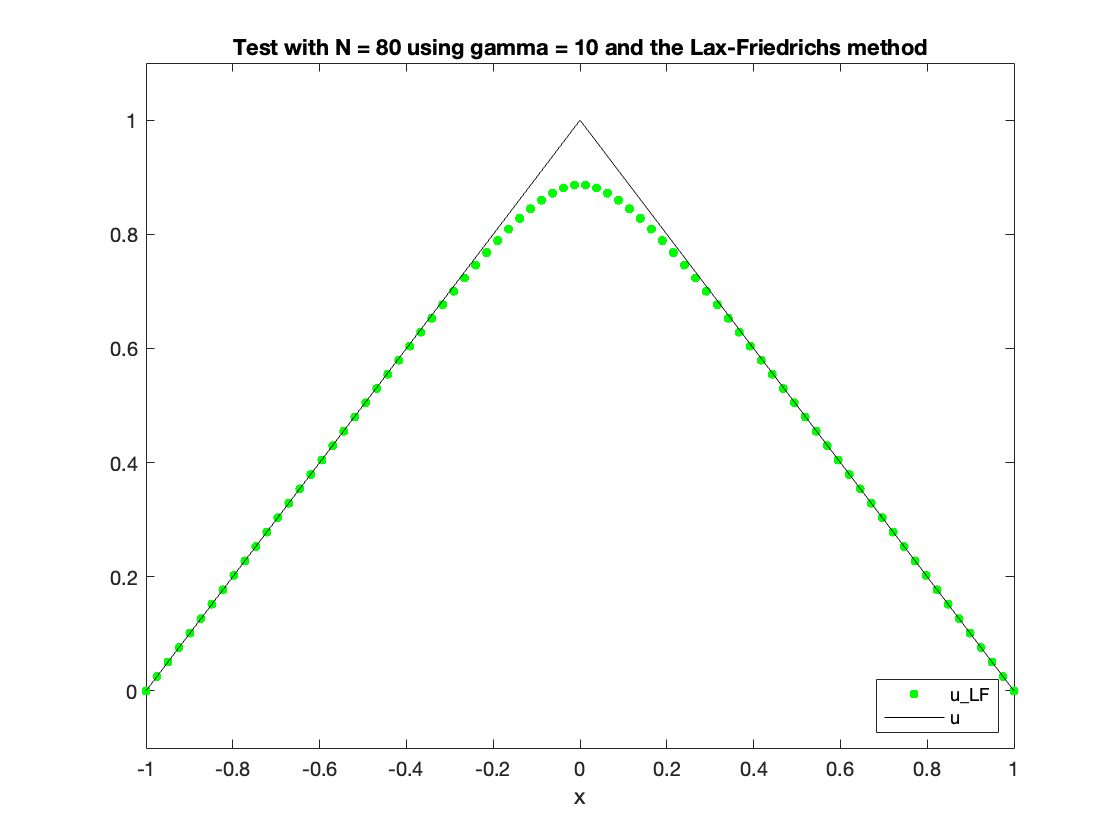} 
\quad 
\includegraphics[width=0.29\textwidth]{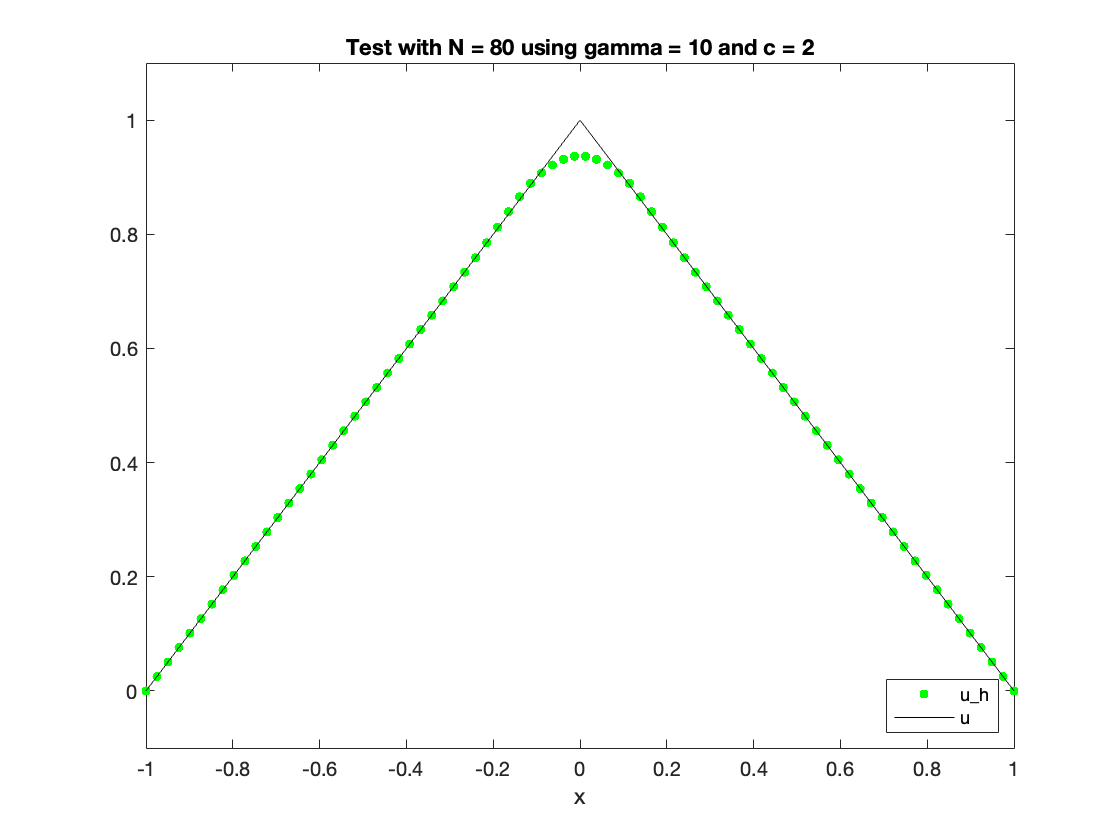} 
\quad 
\includegraphics[width=0.29\textwidth]{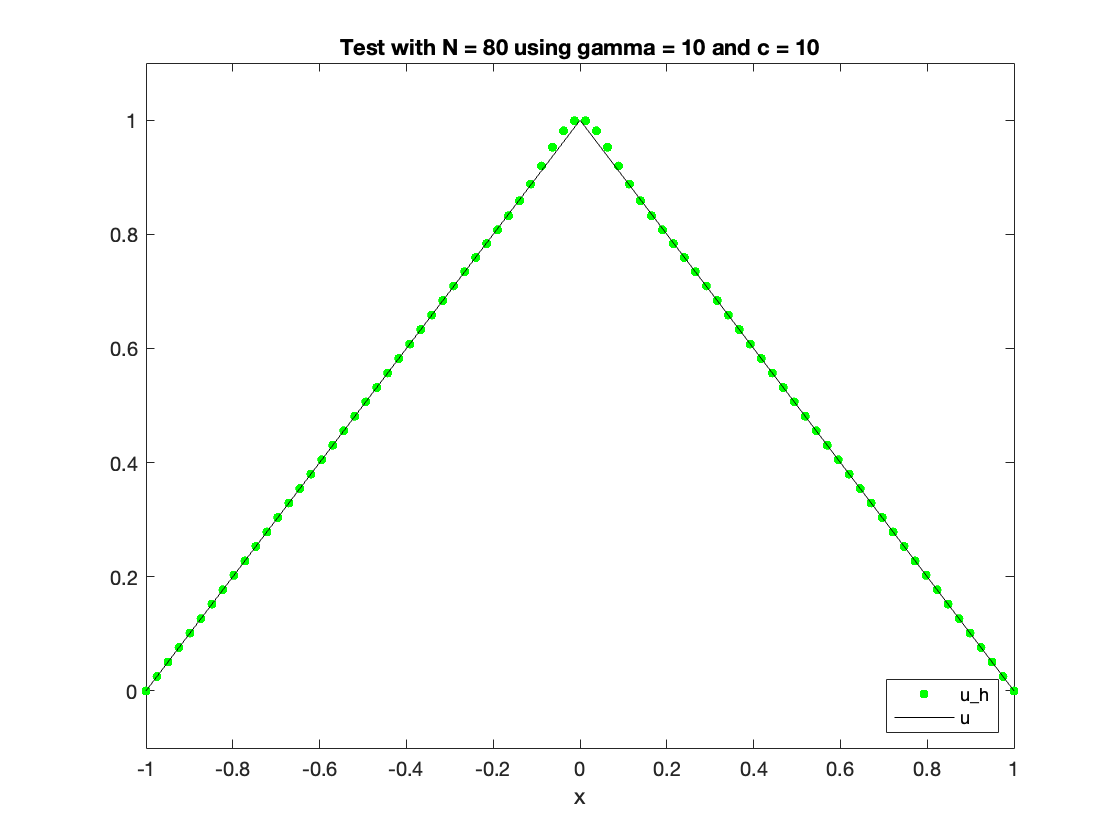} 
\caption{
Approximations for Example 2 in one dimension 
using the Lax-Friedrich's method on the left, 
the modified high-order corrector scheme \eqref{FD_filtered} with $c=2$ in the middle, 
and the modified high-order corrector scheme \eqref{FD_filtered} with $c=10$ on the right.  
The approximation with $c=10$ agrees with the solution to the proposed 
high order scheme \eqref{FD_scheme}.  
All approximations correspond to $h${\em=2.53e-02}.  
Similar results are observed for both boundary conditions \eqref{bc2a} and \eqref{bc2b}.  
}
\label{1D-2plots}
\end{center}
\end{figure}

Consider the problem 
\[
	H[u] = \left| u_x \right| + u - 1 - |x|  
\]
with boundary data chosen such that the viscosity solution is $u(x) = 1 - |x|$.  
We compare the errors for the Lax-Friedrich's method and the proposed high-order method \eqref{FD_scheme} 
in Table~\ref{1D-2a}. 
While both methods appear to converge at a rate of 1, the proposed method is much more accurate 
for the more coarse meshes.  
The results for the modified high-order corrector method with various choices for $c$ can be found in Table~\ref{1D-2b}.
We see in Figure~\ref{1D-2plots} that the approximations become less diffusive as the value for 
$c$ increases.  
Due to the lack of monotonicity, we also see that $c>0$ produces inflection points while the 
Lax-Friedrich's approximation is always concave down consistent with the lower-regularity viscosity solution.  
Consequently, the Lax-Friedrich's method yields a better qualitative solution 
while the proposed non-monotone method yields a more accurate solution near the corner.

\begin{table}[htb] 
{\scriptsize 
\begin{center}
\begin{tabular}{| c | c | c | c | c |}
		\hline
	& \multicolumn{2}{|c|}{Lax-Friedrich's} & \multicolumn{2}{|c|}{Scheme \eqref{FD_scheme}} \\ 
		\hline
	 $h$ & Error & Order & Error & Order \\ 
		\hline
	2.22e-01 & 4.52e-01 &  & 1.27e-01 &  \\ 
		\hline
	1.05e-01 & 3.13e-01 & 0.49 & 6.76e-02 & 0.84 \\ 
		\hline
	5.13e-02 & 1.85e-01 & 0.73 & 3.68e-02 & 0.85 \\ 
		\hline
	2.53e-02 & 1.01e-01 & 0.86 & 1.91e-02 & 0.93 \\ 
		\hline
	1.26e-02 & 5.31e-02 & 0.92 & 9.75e-03 & 0.96 \\ 
		\hline
	6.69e-03 & 2.91e-02 & 0.95 & 5.25e-03 & 0.98 \\ 
		\hline
	4.01e-03 & 1.77e-02 & 0.97 & 3.16e-03 & 0.99 \\ 
		\hline
\end{tabular}
\end{center}
}
\caption{
Approximations for Example 2 in one dimension for the Lax-Friedrich's method and \eqref{FD_scheme} with boundary condition \eqref{bc2a}.  
Similar results hold for the boundary condition \eqref{bc2b}.
}
\label{1D-2a}
\end{table}

\begin{table}[htb] 
{\scriptsize 
\begin{center}
\begin{tabular}{| c | c | c | c | c | c | c | c | c | c |}
		\hline
	& \multicolumn{3}{|c|}{c = 1} & \multicolumn{3}{|c|}{c = 2} & \multicolumn{3}{|c|}{c = 4} \\ 
		\hline
	 $h$ & Cutoff & Error & Order & Cutoff & Error & Order & Cutoff & Error & Order \\ 
		\hline
	2.22e-01 & yes & 2.49e-01 &  & yes & 9.20e-02 &  & no & 1.35e-01 &  \\ 
		\hline
	1.05e-01 & yes & 2.11e-01 & 0.22 & yes & 1.13e-01 & -0.27 & no & 6.73e-02 & 0.93 \\ 
		\hline
	5.13e-02 & yes & 1.35e-01 & 0.62 & yes & 8.48e-02 & 0.39 & yes & 2.58e-02 & 1.33 \\ 
		\hline
	2.53e-02 & yes & 7.59e-02 & 0.81 & yes & 5.10e-02 & 0.72 & yes & 1.16e-02 & 1.13 \\ 
		\hline
	1.26e-02 & yes & 4.06e-02 & 0.90 & yes & 2.81e-02 & 0.85 & yes & 5.40e-03 & 1.09 \\ 
		\hline
	6.69e-03 & yes & 2.24e-02 & 0.94 & yes & 1.57e-02 & 0.92 & yes & 3.46e-03 & 0.70 \\ 
		\hline
	4.01e-03 & yes & 1.36e-02 & 0.97 & yes & 9.66e-03 & 0.95 & yes & 2.26e-03 & 0.83 \\ 
		\hline
\end{tabular}
\end{center}
}
\caption{
Approximations for Example 2 in one dimension for \eqref{FD_filtered} with boundary condition \eqref{bc2b} 
using various smaller values for $c$ that ensure the cutoffs are applied.  
Similar results hold for the boundary condition \eqref{bc2a}.
}
\label{1D-2b}
\end{table}

\subsection{Two-dimensional tests} \label{2D_numerics_sec}

We now consider a series of numerical experiments in two dimensions to test the accuracy of the proposed schemes. 
Overall, the Lax-Friedrich's method has a rate of convergence approaching 1 (asymptotically) in the $\ell^\infty$-norm as expected, 
and it is less accurate than the proposed non-monotone FD methods.  
The proposed high-order method with boundary condition \eqref{bc2a} appears to over-perform and the 
proposed high-order method with boundary condition \eqref{bc2b} appears to under-perform with 
regards to the consistency analysis and observed rates of convergence; 
however, we do see that the boundary condition \eqref{bc2b} typically yields a more accurate approximation 
than the boundary condition \eqref{bc2a} as expected.  
All of the examples satisfy the analytic assumptions made for the differential operator $H$ 
and boundary data in \eqref{HJ}.  
We again first benchmark the methods with a linear problem.  


\subsubsection{Example 1:  linear operator with a smooth solution} 

We first benchmark the scheme by considering the linear problem 
\begin{align*}
	H[u] = u_x + u_y + u - f 
\end{align*}
with $f$ and $g$ chosen such that the exact solution is $u(x,y)=e^{xy}$. 
The results for the Lax-Friedrich's method and the proposed high-order methods can be found in 
Table~\ref{2D-1a} where we see the proposed methods are significantly more accurate.  
The rates when using the boundary condition \eqref{bc2b} are not clear, 
but the method is the most accurate.  
The results for the modified high-order corrector method with various choices for $c$ can be found in Table~\ref{2D-1b} 
where we see a transition region with regards to the impact of the chosen cutoff operators.

\begin{table}[htb] 
{\scriptsize
\begin{center}
\begin{tabular}{| c | c | c | c | c | c | c |}
		\hline
	$h$ & \multicolumn{2}{|c|}{Lax-Friedrich's} & \multicolumn{2}{|c|}{Scheme \eqref{FD_scheme}, BC \eqref{bc2a}} & \multicolumn{2}{|c|}{Scheme \eqref{FD_scheme}, BC \eqref{bc2b}} \\ 
		\hline
	3.14e-01 & 7.95e-02 &  & 3.13e-02 &  & 9.04e-03 &  \\ 
		\hline
	1.49e-01 & 7.19e-02 & 0.14 & 1.53e-02 & 0.96 & 2.77e-03 & 1.58 \\
		\hline
	7.25e-02 & 5.66e-02 & 0.33 & 5.92e-03 & 1.32 & 5.90e-04 & 2.15 \\ 
		\hline
	4.79e-02 & 4.65e-02 & 0.47 & 3.13e-03 & 1.54 & 2.20e-04 & 2.38 \\ 
		\hline
	3.58e-02 & 3.95e-02 & 0.56 & 1.92e-03 & 1.68 & 1.10e-04 & 2.37 \\ 
		\hline
	2.86e-02 & 3.46e-02 & 0.59 & 1.29e-03 & 1.75 & 8.77e-05 & 1.01 \\ 
		\hline
	2.38e-02 & 3.06e-02 & 0.66 & 9.29e-04 & 1.79 & 7.30e-05 & 1.00 \\ 
		\hline
	1.78e-02 & 2.52e-02 & 0.68 & 5.45e-04 & 1.84 & 5.09e-05 & 1.24 \\ 
		\hline
	1.42e-02 & 2.13e-02 & 0.73 & 3.58e-04 & 1.87 & 3.70e-05 & 1.43 \\ 
		\hline
\end{tabular}
\end{center}
}
\caption{
Approximations for Example 1 in two dimensions for the Lax-Friedrich's method and the proposed high-order scheme \eqref{FD_scheme} 
with the various choices for the auxiliary boundary condition.  
}
\label{2D-1a}
\end{table}

\begin{table}[htb] 
{\scriptsize 
\begin{center}
\begin{tabular}{| c | c | c | c | c | c | c | c | c | c |}
		\hline
	& \multicolumn{3}{|c|}{c = 0.1} & \multicolumn{3}{|c|}{c = 1} & \multicolumn{3}{|c|}{c = 1.2} \\ 
		\hline
	 $h$ & Cutoff & Error & Order & Cutoff & Error & Order & Cutoff & Error & Order \\ 
		\hline
	3.14e-01 & yes & 7.83e-02 &  & no & 9.04e-03 &  & no & 9.04e-03 &  \\ 
		\hline
	1.49e-01 & yes & 7.03e-02 & 0.14 & no & 2.77e-03 & 1.58 & no & 2.77e-03 & 1.58 \\ 
		\hline
	7.25e-02 & yes & 5.48e-02 & 0.35 & no & 5.90e-04 & 2.15 & no & 5.90e-04 & 2.15 \\ 
		\hline
	4.79e-02 & yes & 4.49e-02 & 0.48 & yes & 2.20e-04 & 2.38 & no & 2.20e-04 & 2.38 \\ 
		\hline
	3.58e-02 & yes & 3.84e-02 & 0.54 & yes & 1.03e-02 & -13.19 & no & 1.10e-04 & 2.37 \\ 
		\hline
	2.86e-02 & yes & 3.33e-02 & 0.63 & yes & 1.97e-02 & -2.85 & yes & 1.15e-03 & -10.40 \\ 
		\hline
	2.38e-02 & yes & 2.98e-02 & 0.61 & yes & 2.36e-02 & -0.99 & yes & 6.03e-03 & -9.00 \\ 
		\hline
	1.78e-02 & yes & 2.44e-02 & 0.69 & yes & 2.15e-02 & 0.32 & yes & 1.39e-02 & -2.87 \\ 
		\hline
	1.42e-02 & yes & 2.07e-02 & 0.74 & yes & 1.86e-02 & 0.66 & yes & 1.60e-02 & -0.64 \\ 
		\hline
\end{tabular}
\end{center}
}
\caption{
Approximations for Example 1 in two dimensions for \eqref{FD_filtered} with boundary condition \eqref{bc2b} 
using various smaller values for $c$ that ensure the cutoffs are applied.  
Similar results hold for the boundary condition \eqref{bc2a}.
}
\label{2D-1b}
\end{table}


\subsubsection{Example 2:  nonlinear $C^1 \backslash C^2$ operator with a smooth solution} 

We next consider the nonlinear problem 
\begin{align*}
	H[u] = \sqrt{u_x^2 + u_y^2} + u - f 
\end{align*}
with $f$ and $g$ chosen such that the exact solution is $u(x,y)=e^{xy}$.
The results for the Lax-Friedrich's method and the proposed high-order method can be found in 
Table~\ref{2D-2a} where we see the proposed method is more accurate.  
The Lax-Friedrich's method is also in the pre-asymptotic region for measuring 
the rates of convergence.  
The results for the modified high-order corrector method with various choices for $c$ can be found in Table~\ref{2D-2b} 
where we see a transition region with regards to the impact of the chosen cutoff operators.  
The slow convergence rate of the Lax-Friedrich's method has a major impact on the accuracy 
of the modified method for finer meshes.  

\begin{table}[htb] 
{\scriptsize
\begin{center}
\begin{tabular}{| c | c | c | c | c | c | c |}
		\hline
	$h$ & \multicolumn{2}{|c|}{Lax-Friedrich's} & \multicolumn{2}{|c|}{Scheme \eqref{FD_scheme}, BC \eqref{bc2a}} & \multicolumn{2}{|c|}{Scheme \eqref{FD_scheme}, BC \eqref{bc2b}} \\ 
		\hline
	3.14e-01 & 7.29e-02 &  & 2.74e-02 &  & 5.84e-03 &  \\ 
		\hline
	1.49e-01 & 6.57e-02 & 0.14 & 1.37e-02 & 0.93 & 1.45e-03 & 1.87 \\ 
		\hline
	7.25e-02 & 5.28e-02 & 0.31 & 6.05e-03 & 1.13 & 4.69e-04 & 1.57 \\ 
		\hline
	4.79e-02 & 4.51e-02 & 0.38 & 3.30e-03 & 1.46 & 2.57e-04 & 1.45 \\ 
		\hline
	3.58e-02 & 3.95e-02 & 0.46 & 2.05e-03 & 1.63 & 1.70e-04 & 1.41 \\ 
		\hline
	2.86e-02 & 3.52e-02 & 0.50 & 1.39e-03 & 1.71 & 1.32e-04 & 1.12 \\ 
		\hline
	2.38e-02 & 3.19e-02 & 0.54 & 1.01e-03 & 1.77 & 1.05e-04 & 1.25 \\ 
		\hline
	1.78e-02 & 2.68e-02 & 0.60 & 5.95e-04 & 1.82 & 7.06e-05 & 1.37 \\ 
		\hline
	1.42e-02 & 2.32e-02 & 0.65 & 3.92e-04 & 1.86 & 5.06e-05 & 1.49 \\ 
		\hline
\end{tabular}
\end{center}
}
\caption{
Approximations for Example 2 in two dimensions for the Lax-Friedrich's method and the proposed high-order scheme \eqref{FD_scheme} 
with the various choices for the auxiliary boundary condition.  
}
\label{2D-2a}
\end{table}

\begin{table}[htb] 
{\scriptsize 
\begin{center}
\begin{tabular}{| c | c | c | c | c | c | c | c | c | c |}
		\hline
	& \multicolumn{3}{|c|}{c = 0.1, BC \eqref{bc2a}} & \multicolumn{3}{|c|}{c = 1, BC \eqref{bc2a}} & \multicolumn{3}{|c|}{c = 1, BC \eqref{bc2b}} \\ 
		\hline
	 $h$ & Cutoff & Error & Order & Cutoff & Error & Order & Cutoff & Error & Order \\ 
		\hline
	3.14e-01 & yes & 5.90e-02 &  & no & 2.74e-02 &  & no & 5.84e-03 &  \\ 
		\hline
	1.49e-01 & yes & 5.77e-02 & 0.03 & no & 1.37e-02 & 0.93 & no & 1.45e-03 & 1.87 \\ 
		\hline
	7.25e-02 & yes & 4.98e-02 & 0.20 & no & 6.05e-03 & 1.13 & no & 4.69e-04 & 1.57 \\ 
		\hline
	4.79e-02 & yes & 4.31e-02 & 0.35 & no & 3.30e-03 & 1.46 & no & 2.57e-04 & 1.45 \\ 
		\hline
	3.58e-02 & yes & 3.77e-02 & 0.45 & yes & 5.85e-03 & -1.96 & yes & 7.78e-03 & -11.68 \\ 
		\hline
	2.86e-02 & yes & 3.35e-02 & 0.52 & yes & 9.23e-03 & -2.02 & yes & 1.68e-02 & -3.41 \\ 
		\hline
	2.38e-02 & yes & 3.04e-02 & 0.52 & yes & 1.14e-02 & -1.14 & yes & 2.29e-02 & -1.67 \\ 
		\hline
	1.78e-02 & yes & 2.56e-02 & 0.60 & yes & 1.33e-02 & -0.54 & yes & 2.06e-02 & 0.36 \\ 
		\hline
	1.42e-02 & yes & 2.22e-02 & 0.63 & yes & 1.36e-02 & -0.09 & yes & 1.82e-02 & 0.54 \\ 
		\hline
\end{tabular}
\end{center}
}
\caption{
Approximations for Example 2 in two dimensions for \eqref{FD_filtered} with various boundary conditions 
using various smaller values for $c$ that ensure the cutoffs are applied.    
Similar results hold for the boundary condition \eqref{bc2b} when $c=0.1$.
}
\label{2D-2b}
\end{table}


\subsubsection{Example 3: nonlinear Lipschitz operator with a smooth solution} 

We next consider the nonlinear problem 
\begin{align*}
	H[u] = |u_x| + |u_y| + |u| + 2u - f
\end{align*}
with $f$ and $g$ chosen such that the exact solution is $u(x,y)=\cos{(\pi x)}\cos{(\pi y)} - 0.5$.
The results for the Lax-Friedrich's method and the proposed high-order method can be found in 
Table~\ref{2D-3a} where we see the proposed method is more accurate.  
The results for the modified high-order corrector method with various choices for $c$ can be found in Table~\ref{2D-3b} 
where we see a transition region with regards to the impact of the chosen cutoff operators.

\begin{table}[htb] 
{\scriptsize
\begin{center}
\begin{tabular}{| c | c | c | c | c | c | c |}
		\hline
	$h$ & \multicolumn{2}{|c|}{Lax-Friedrich's} & \multicolumn{2}{|c|}{Scheme \eqref{FD_scheme}, BC \eqref{bc2a}} & \multicolumn{2}{|c|}{Scheme \eqref{FD_scheme}, BC \eqref{bc2b}} \\ 
		\hline
	3.14e-01 & 4.48e-01 &  & 4.84e-01 &  & 5.25e-01 &  \\ 
		\hline
	1.49e-01 & 4.08e-01 & 0.13 & 1.69e-01 & 1.41 & 9.46e-02 & 2.29 \\ 
		\hline
	7.25e-02 & 3.44e-01 & 0.23 & 3.13e-02 & 2.35 & 1.21e-02 & 2.86 \\ 
		\hline
	4.79e-02 & 2.80e-01 & 0.50 & 1.24e-02 & 2.24 & 4.10e-03 & 2.62 \\ 
		\hline
	3.58e-02 & 2.32e-01 & 0.65 & 6.68e-03 & 2.12 & 2.68e-03 & 1.45 \\ 
		\hline
	2.86e-02 & 1.97e-01 & 0.74 & 4.20e-03 & 2.06 & 1.88e-03 & 1.57 \\ 
		\hline
	2.38e-02 & 1.70e-01 & 0.79 & 2.90e-03 & 2.02 & 1.39e-03 & 1.64 \\ 
		\hline
	1.78e-02 & 1.33e-01 & 0.84 & 1.62e-03 & 2.00 & 8.49e-04 & 1.71 \\ 
		\hline
	1.42e-02 & 1.09e-01 & 0.88 & 1.04e-03 & 1.99 & 5.70e-04 & 1.78 \\ 
		\hline
\end{tabular}
\end{center}
}
\caption{
Approximations for Example 3 in two dimensions for the Lax-Friedrich's method and the proposed high-order scheme \eqref{FD_scheme} 
with the various choices for the auxiliary boundary condition.  
}
\label{2D-3a}
\end{table}

\begin{table}[htb] 
{\scriptsize 
\begin{center}
\begin{tabular}{| c | c | c | c | c | c | c | c | c | c |}
		\hline
	& \multicolumn{3}{|c|}{c = 1} & \multicolumn{3}{|c|}{c = 4} & \multicolumn{3}{|c|}{c = 6} \\ 
		\hline
	 $h$ & Cutoff & Error & Order & Cutoff & Error & Order & Cutoff & Error & Order \\ 
		\hline
	3.14e-01 & yes & 4.05e-01 &  & no & 4.84e-01 &  & no & 4.84e-01 &  \\ 
		\hline
	1.49e-01 & yes & 3.12e-01 & 0.35 & no & 1.69e-01 & 1.41 & no & 1.69e-01 & 1.41 \\ 
		\hline
	7.25e-02 & yes & 2.76e-01 & 0.17 & yes & 6.59e-02 & 1.31 & no & 3.13e-02 & 2.35 \\ 
		\hline
	4.79e-02 & yes & 2.33e-01 & 0.41 & yes & 9.10e-02 & -0.78 & no & 1.24e-02 & 2.24 \\ 
		\hline
	3.58e-02 & yes & 1.96e-01 & 0.58 & yes & 8.98e-02 & 0.05 & yes & 1.93e-02 & -1.51 \\ 
		\hline
	2.86e-02 & yes & 1.68e-01 & 0.69 & yes & 8.28e-02 & 0.36 & yes & 2.61e-02 & -1.34 \\ 
		\hline
	2.38e-02 & yes & 1.46e-01 & 0.76 & yes & 7.52e-02 & 0.52 & yes & 2.79e-02 & -0.35 \\ 
		\hline
	1.78e-02 & yes & 1.15e-01 & 0.82 & yes & 6.22e-02 & 0.65 & yes & 2.68e-02 & 0.14 \\ 
		\hline
	1.42e-02 & yes & 9.51e-02 & 0.86 & yes & 5.25e-02 & 0.75 & yes & 2.42e-02 & 0.45 \\ 
		\hline
\end{tabular}
\end{center}
}
\caption{
Approximations for Example 3 in two dimensions for \eqref{FD_filtered} with boundary condition \eqref{bc2a} 
using various smaller values for $c$ that ensure the cutoffs are applied.  
Similar results hold for the boundary condition \eqref{bc2b}.
}
\label{2D-3b}
\end{table} 


\subsubsection{Example 4:  nonlinear Lipschitz operator with a non-smooth solution} 
We lastly consider the nonlinear problem 
\begin{align*}
	H[u] = |u_x| + 2u_x  + u - f 
\end{align*}
with $f$ and $g$ chosen such that the exact solution is $u(x,y) = | x - 0.2 |$. 
The results for the Lax-Friedrich's method and the proposed high-order method can be found in 
Table~\ref{2D-4a} where we see the proposed method is more accurate.  
The results for the modified high-order corrector method with various choices for $c$ can be found in Table~\ref{2D-4b} 
where we see a transition region with regards to the impact of the chosen cutoff operators.  
As expected, the methods are only first-order accurate due to the lower regularity of the viscosity solution.  

\begin{table}[htb] 
{\scriptsize
\begin{center}
\begin{tabular}{| c | c | c | c | c | c | c |}
		\hline
	$h$ & \multicolumn{2}{|c|}{Lax-Friedrich's} & \multicolumn{2}{|c|}{Scheme \eqref{FD_scheme}, BC \eqref{bc2a}} & \multicolumn{2}{|c|}{Scheme \eqref{FD_scheme}, BC \eqref{bc2b}} \\ 
		\hline
	3.14e-01 & 4.76e-01 &  & 1.46e-01 &  & 1.36e-01 &  \\ 
		\hline
	1.49e-01 & 3.11e-01 & 0.57 & 6.58e-02 & 1.07 & 6.91e-02 & 0.90 \\ 
		\hline
	7.25e-02 & 1.72e-01 & 0.82 & 3.48e-02 & 0.89 & 3.63e-02 & 0.89 \\ 
		\hline
	4.79e-02 & 1.18e-01 & 0.91 & 2.36e-02 & 0.94 & 2.47e-02 & 0.94 \\ 
		\hline
	3.58e-02 & 9.02e-02 & 0.93 & 1.79e-02 & 0.96 & 1.87e-02 & 0.96 \\ 
		\hline
	2.86e-02 & 7.29e-02 & 0.95 & 1.44e-02 & 0.97 & 1.50e-02 & 0.97 \\ 
		\hline
	2.38e-02 & 6.11e-02 & 0.95 & 1.20e-02 & 0.97 & 1.26e-02 & 0.97 \\ 
		\hline
	1.78e-02 & 4.62e-02 & 0.96 & 9.05e-03 & 0.98 & 9.46e-03 & 0.98 \\ 
		\hline
	1.42e-02 & 3.72e-02 & 0.97 & 7.26e-03 & 0.98 & 7.58e-03 & 0.98 \\ 
		\hline
\end{tabular}
\end{center}
}
\caption{
Approximations for Example 4 in two dimensions for the Lax-Friedrich's method and the proposed high-order scheme \eqref{FD_scheme} 
with the various choices for the auxiliary boundary condition.  
}
\label{2D-4a}
\end{table}

\begin{table}[htb] 
{\scriptsize 
\begin{center}
\begin{tabular}{| c | c | c | c | c | c | c | c | c | c |}
		\hline
	& \multicolumn{3}{|c|}{c = 0.1} & \multicolumn{3}{|c|}{c = 1} & \multicolumn{3}{|c|}{c = 3} \\ 
		\hline
	 $h$ & Cutoff & Error & Order & Cutoff & Error & Order & Cutoff & Error & Order \\ 
		\hline
	3.14e-01 & yes & 4.85e-01 &  & yes & 3.07e-01 &  & no & 1.46e-01 &  \\ 
		\hline
	1.49e-01 & yes & 2.98e-01 & 0.65 & yes & 1.94e-01 & 0.62 & no & 6.58e-02 & 1.07 \\ 
		\hline
	7.25e-02 & yes & 1.65e-01 & 0.82 & yes & 1.09e-01 & 0.80 & no & 3.48e-02 & 0.89 \\ 
		\hline
	4.79e-02 & yes & 1.14e-01 & 0.90 & yes & 7.48e-02 & 0.90 & no & 2.36e-02 & 0.94 \\ 
		\hline
	3.58e-02 & yes & 8.67e-02 & 0.93 & yes & 5.71e-02 & 0.93 & no & 1.79e-02 & 0.96 \\ 
		\hline
	2.86e-02 & yes & 7.00e-02 & 0.94 & yes & 4.62e-02 & 0.94 & yes & 1.44e-02 & 0.97 \\ 
		\hline
	2.38e-02 & yes & 5.88e-02 & 0.95 & yes & 3.87e-02 & 0.95 & yes & 1.20e-02 & 0.97 \\ 
		\hline
	1.78e-02 & yes & 4.45e-02 & 0.96 & yes & 2.93e-02 & 0.96 & yes & 9.05e-03 & 0.98 \\ 
		\hline
	1.42e-02 & yes & 3.58e-02 & 0.97 & yes & 2.36e-02 & 0.97 & yes & 7.26e-03 & 0.98 \\ 
		\hline
\end{tabular}
\end{center}
}
\caption{
Approximations for Example 4 in two dimensions for \eqref{FD_filtered} with boundary condition \eqref{bc2a} 
using various smaller values for $c$ that ensure the cutoffs are applied.  
Similar results hold for the boundary condition \eqref{bc2b}.
}
\label{2D-4b}
\end{table}



\section{Conclusion} \label{conclusion_sec}

In this paper we have proposed a new high-order approximation method for non-degenerate stationary 
Hamilton-Jacobi problems with Dirichlet boundary data as well as a modified version for which we can 
guarantee admissibility and stability.   
The methods both appear to be more accurate than the Lax-Friedrich's method and capable of 
accuracy beyond the Godunov barrier for monotone methods.  
The main idea is to use a numerical moment which is less diffusive than the typical numerical viscosity approach.  
Rewriting the numerical moment, we see that the new methods can be regarded as high-order corrections 
to the Lax-Friedrich's method.  
Using results for Lax-Friedrich's methods, we can prove convergence of the proposed modified high-order corrector method.  
In the numerical tests, the modified high-order corrector scheme often agreed with the proposed unmodified version ensuring 
the non-monotone method was converging.  
It is easy to check whether or not the cutoffs used to modify the high-order corrector method were applied.  
Future work will seek to directly analyze the high-order formulation without having to rely upon 
the cutoff operators for proving various admissibility, stability, and convergence results.  
The non-monotone methods performed well in all numerical experiments offering higher accuracy than the 
corresponding Lax-Friedrich's method.  
Unfortunately, the lack of monotonicity creates major analytic hurdles making the modified scheme 
and its underlying analysis a strong intermediate step as we seek to develop the new analytic techniques 
needed to analyze more general non-monotone methods.  


\bibliographystyle{elsarticle-num}

\end{document}